\documentclass[11pt]{amsart}

\usepackage[T1]{fontenc}
\usepackage[utf8]{inputenc}
\usepackage{amsmath,amssymb,amsfonts,amsthm,mathrsfs,mathtools}
\usepackage[expansion=false]{microtype}
\usepackage{xcolor}
\usepackage[colorlinks=true,citecolor=blue]{hyperref}

\allowdisplaybreaks[1]
\newtheoremstyle{fancy}{}{}{\itshape}{}{\textbf\bgroup}{.\egroup}{ }{}
\newtheoremstyle{fanci}{}{}{\rmfamily}{}{\textbf\bgroup}{.\egroup}{ }{}
\numberwithin{equation}{section}
\swapnumbers
\theoremstyle{fancy}
\newtheorem{theorem}[equation]{Theorem}
\newtheorem{proposition}[equation]{Proposition}
\newtheorem{lemma}[equation]{Lemma}
\newtheorem{corollary}[equation]{Corollary}
\theoremstyle{fanci}
\newtheorem{definition}[equation]{Definition}
\newtheorem{example}[equation]{Example}

\theoremstyle{fancy}
\newtheorem*{introindex}{Theorem 1}
\newtheorem*{introflag}{Theorem 2}
\newtheorem*{introsphere}{Theorem 3}

\DeclareMathOperator{\Proj}{Proj}
\newcommand{\Ent}{\operatorname{Ent}}
\newcommand{\Hess}{\operatorname{Hess}}
\newcommand{\Ric}{\operatorname{Ric}}
\newcommand{\Hom}{\operatorname{Hom}}

\newcommand{\Fl}{\operatorname{Fl}}

\title[Spectral instability of Wasserstein entropy]{Spectral Instability of Wasserstein Entropy on Homogeneous Spaces}

\author[G. C. Brito]{Geovane C. Brito}
\address{Geovane C. Brito, Institute of Mathematics, Department of Applied Mathematics, Universidade Estadual de Campinas, 13.083-859 Campinas - SP, Brazil}
\email{gcbritomath@gmail.com}
\author[M. Y. Maruyama]{Matheus Y. Maruyama}
\address{Matheus Y. Maruyama, Institute of Mathematics, Department of Applied Mathematics, Universidade Estadual de Campinas, 13.083-859 Campinas - SP, Brazil}
\email{maruyama@ime.unicamp.br}
\author[C. S. Rodrigues]{Christian S. Rodrigues}
\address{Christian S. Rodrigues, Institute of Mathematics, Department of Applied Mathematics, Universidade Estadual de Campinas, 13.083-859 Campinas - SP, Brazil, and Max-Planck-Institute for Mathematics in the Sciences, Inselstr. 22, 04103 Leipzig, Germany}
\email{rodrigues@ime.unicamp.br}

\date{}
\hypersetup{
 pdftitle={Spectral Instability of Wasserstein Entropy on Homogeneous Spaces},
 pdfauthor={Geovane C. Brito, Matheus Y. Maruyama, Christian S. Rodrigues}
}

\begin{document}

\begin{abstract}
Let $M$ be a compact homogeneous Riemannian manifold, and let $m$ be its invariant probability measure. At $m$, the Wasserstein entropy Hessian has a sharp lower bound given by $\lambda_1$, the first positive eigenvalue of the Laplace–Beltrami operator on $M$. We study how this bound changes under perturbations of the density by counting the generalized Hessian eigenvalues below $\lambda_1$. Using a Schur complement and Peter–Weyl triple products, we reduce this count to the first eigenspace and its interaction with higher modes, including isotropy multiplicities.

We apply this method to two important families of homogeneous manifolds. On the unit round sphere $S^d$, $d\geq 2$, we prove that the bound $\lambda_1=d$ holds only at the uniform density and obtain an exact local index formula for every nonzero smooth mean-zero perturbation $\eta$. Writing its degree-two component as $x^{\mathsf T}Qx$, the index at $1+\varepsilon\eta$ is $n_-(\varepsilon Q)+\dim\ker Q$ for sufficiently small nonzero $\varepsilon$. On the normal complete flags $SU(N)/T^{N-1}$, we find a perturbation with complementary indices $N(N-2)^2+2N-1$ and $(N-2)(3N-1)$ under opposite signs, whose proportions tend to one and zero as $N$ grows. A resonant example on $S^1\times S^2$ shows that spherical rigidity does not hold on every homogeneous manifold.
\end{abstract}

\maketitle

\section{Introduction}\label{sec:introduction}

At the uniform probability measure of a compact homogeneous manifold, the Wasserstein Hessian of entropy has a sharp lower bound determined by the first positive Laplace eigenvalue, and this bound is attained precisely in the directions given by gradients of first eigenfunctions. This provides a distinguished reference configuration in which both the optimal constant and its equality directions are explicitly known.

What happens to this comparison when mass is redistributed? A change of density alters both the Otto metric and the Wasserstein entropy Hessian, allowing spectral modes that were previously orthogonal to interact. The homogeneous bound may then fail in some directions, while still holding in others. Our question is when this bound still holds and, when it fails, how many independent directions violate it.

To answer this question, let $F=G/H$ be a compact connected homogeneous Riemannian manifold of positive dimension, with normalized invariant measure $m$, and let $\rho m$ be a probability measure with smooth positive density $\rho$.

For smooth potentials $u,v$, we define the Otto metric and the Wasserstein entropy Hessian by
\begin{equation}\label{eq:intro-forms}
\begin{aligned}
\mathfrak g_\rho(u,v)
 &=\int_F\rho\langle\nabla u,\nabla v\rangle\,dm,\\
\mathfrak h_\rho(u,v)
 &=\int_F\rho\bigl(\langle\Hess u,\Hess v\rangle
       +\Ric(\nabla u,\nabla v)\bigr)\,dm.
\end{aligned}
\end{equation}

The smooth Otto calculus and the second-variation formulas for entropy were studied in \cite{Otto2001,Lott2008}, while their connection with Ricci curvature bounds and displacement convexity was developed in \cite{Bakry1991,RenesseSturm2005,LottVillani2009}. In this paper, we focus on a different spectral property at a particular measure, instead of a global lower bound for curvature.

Let $L=-\Delta$ be the Laplace–Beltrami operator on $F$, and let $\lambda_1>0$ be its first positive eigenvalue. The integrated Bochner identity gives
\begin{equation}
  \label{eq:intro-gap}
  \mathfrak{h}_{1}(u,u)\geq \lambda_{1} \mathfrak{g}_{1}(u,u),
\end{equation}
with equality precisely for functions in $E_{1}$, the eigenspace associated with $\lambda_{1}$, modulo constants. Thus, $\lambda_1$ gives a natural spectral threshold, determined by the geometry of the manifold, since it is the smallest generalized eigenvalue of the entropy Hessian at $m$. For $\rho\neq 1$, this threshold remains fixed, allowing us to study whether changing the density causes the same optimal inequality to fail, without changing the Laplacian itself. In particular, the negative
index of
\begin{equation}\label{eq:intro-defect}
\mathfrak d_\rho=\mathfrak h_\rho-\lambda_1\mathfrak g_\rho
\end{equation}
counts, with multiplicity, the generalized eigenvalues of the entropy Hessian that are below $\lambda_1$.

Throughout the paper, we use the term \emph{spectral instability} to describe eigenvalues crossing the fixed positive threshold $\lambda_1$. For example, on the unit round sphere, $\Ric=(d-1)g$ implies that $\mathfrak h_\rho\geq(d-1)\mathfrak g_\rho$ for every smooth positive density. Therefore, the sharper bound $d$ may fail even though the entropy Hessian remains strictly positive.

The homogeneous structure allows us to compute this index using representation theory, since Peter--Weyl theory decomposes functions on $G/H$ into irreducible representations together with their $H$-fixed multiplicity spaces. When the density changes, the spectral blocks interact through products of eigenfunctions, and we use the Clebsch--Gordan decomposition to identify the possible interactions, while the corresponding triple integrals and Laplace eigenvalues determine which of them contribute to the defect. This approach is based on classical harmonic analysis \cite{BrockerDieck1985,Farashahi2017}, with related spectral coordinates in Wasserstein geometry appearing in \cite{GomesRodriguesSanMartin2026,GomesRodriguesInvariant2026}.

To compute the index, we combine these representation-theoretic tools with the methods of Maddocks and Kato \cite{Maddocks1988,Kato1995}, obtaining explicit coupling operators, finite-channel tests, and signatures that account for the multiplicities of the homogeneous space. Although the abstract index reduction holds on every closed connected Riemannian manifold of positive dimension with normalized volume, homogeneity provides the trace identity and the representation coordinates that allow us to evaluate the reduction geometrically.

These ideas bring us to the main question of this paper, which connects the behavior of the entropy Hessian under perturbations with the representation theory of homogeneous spaces.

\medskip

\textit{For which perturbations of an invariant measure does the sharp
Wasserstein entropy-Hessian inequality fail, and what can the
representation theory of a homogeneous space tell us about the number
of directions in which it fails?}

\medskip
\medskip

The answer to this question can vary even among homogeneous spaces. Some nonconstant density perturbations are invisible to the critical spectral row because of resonant couplings, as we show below with an explicit example on $S^1\times S^2$. On round spheres, however, every nonconstant smooth positive density is detected, while on complete flags, changing the sign of a perturbation can lead to very different negative indices. These examples show why it is important to understand not only whether the bound fails, but also the \emph{order} at which instability appears and the \emph{number of independent directions} in which it occurs.

\subsection*{Statement of the main results}
The first result reduces the instability problem to the first positive
Laplace eigenspace. Let $\eta$ be a smooth real function of zero mean
and consider the density $\rho_\varepsilon=1+\varepsilon\eta$.
Because the defect depends linearly on the density, its variation at
$m$ is the form
\[
\mathfrak d_\eta(u,v)
=\left.\frac{d}{d\varepsilon}\right|_{\varepsilon=0}
\mathfrak d_{1+\varepsilon\eta}(u,v).
\]

Since the unperturbed defect vanishes on $E_1$ and is strictly positive on the higher Laplace modes, we can split the real mean-zero Sobolev space into a critical part and its orthogonal complement:
\[
\mathcal V=H^2(F,\mathbb R)\cap L^2_0(F,m)
=E_1\oplus\mathcal V_+,
\qquad \mathcal V_+=\mathcal V\cap E_1^\perp.
\]

The variation defines a symmetric operator $M_\eta$ when restricted to $E_1$, and a mixed operator $B_\eta$ through its pairing with the higher-mode complement. Using the $L^2(m)$ inner product on $E_1$ and duality on $\mathcal V_+$, we characterize these operators by
\begin{equation}\label{eq:intro-operators}
\begin{aligned}
\langle M_\eta u,w\rangle&=\mathfrak d_\eta(u,w),
&&u,w\in E_1,\\
\langle B_\eta u,v\rangle&=\mathfrak d_\eta(u,v),
&&u\in E_1,\quad v\in\mathcal V_+.
\end{aligned}
\end{equation}

For a finite-dimensional symmetric operator $M$, we denote by $n_-(M)$ and $n_+(M)$ the numbers of negative and positive eigenvalues, counted with multiplicity, while the index and nullity of $\mathfrak d_\rho$ are defined on real potentials modulo constants. With this notation, we obtain the following results.
\begin{introindex}
Let $F$ be a closed connected Riemannian manifold of positive dimension, equipped with its normalized Riemannian volume, and let $\eta$ be a smooth real mean-zero perturbation. Then there exists $\varepsilon_0>0$, depending on the metric and on $\eta$, such that, for every $0<|\varepsilon|<\varepsilon_0$ with $\rho_\varepsilon>0$, we have
\begin{equation}\label{eq:intro-index}
\begin{aligned}
\operatorname{ind}_-(\mathfrak d_{\rho_\varepsilon})
&=n_-(\varepsilon M_\eta)
+\operatorname{rank}(B_\eta|_{\ker M_\eta}),\\
\operatorname{nul}(\mathfrak d_{\rho_\varepsilon})
&=\dim(\ker M_\eta\cap\ker B_\eta).
\end{aligned}
\end{equation}
In particular, the negative index is at most $\dim E_1$.
\end{introindex}
The theorem is proved in Section ~\ref{sec:coercive}, where we also establish Proposition~\ref{prop:spectral-gram}, which describes the spectral structure of the index formula. Its first term measures the splitting within the critical eigenspace, while the second counts the negative directions arising from its coupling with higher modes, even when the first-order restriction vanishes. On homogeneous spaces, both terms can be computed using Peter--Weyl triple products and a positive Gram series, and for perturbations with finite Peter--Weyl support, Proposition~\ref{prop:finite-channel-test} reduces the rank calculations to a finite problem that includes the $H$-fixed multiplicities. In this way, the formula distinguishes first-order splitting, second-order instability, and complete resonance of the critical spectral channels.

On complete flag manifolds, the critical eigenspace contains
multiple copies of the adjoint representation, making the
multiplicity structure indispensable to the index calculation.
We consider the normal metric on
$\Fl_N=SU(N)/T^{N-1}$ induced by
$-\operatorname{tr}(XY)$ on $\mathfrak{su}(N)$, and write $m_N$
for the invariant probability measure. For each $N\ge3$, the
following mean-zero perturbation direction admits a closed
calculation:
\begin{equation}\label{eq:intro-flag-direction}
\eta_N(gT)=N^2|g_{NN}|^2-N,
\qquad \rho_{\varepsilon,N}=1+\varepsilon\eta_N.
\end{equation}

\begin{introflag}
For every $N\ge3$, the normal complete flag $\Fl_N$ satisfies
\[
\lambda_1(\Fl_N)=2N,
\qquad \dim E_1=(N^2-1)(N-1).
\]
For each fixed $N$ and all sufficiently small admissible
$\varepsilon\ne0$, the negative index of the entropy-Hessian defect
is
\begin{equation}\label{eq:intro-flag-index}
\operatorname{ind}_-(\mathfrak d_{\rho_{\varepsilon,N}})
=\begin{cases}
N(N-2)^2+2N-1,&\varepsilon>0,\\
(N-2)(3N-1),&\varepsilon<0.
\end{cases}
\end{equation}
\end{introflag}

The two signs give complementary signatures on $E_1$, although their proportions behave very differently as $N$ increases: for one sign, almost every critical direction becomes negative, while for the other, the proportion of negative directions tends to zero. This difference follows from the explicit cubic Haar moment and the tensor-product factorization of the first variation, with the full calculations, including the multiplicity-space eigenvalues, given in Section~\ref{sec:complete-flags}, where Theorem~\ref{thm:flag-index} is proved.

Round spheres show a different behavior, since their harmonic product decomposition is particularly simple, but the critical couplings still detect every nonconstant density. To study this case, consider the unit sphere $S^d=SO(d+1)/SO(d)$, with $d\geq 2$, and let $\mathscr H_n$ denote the space of real spherical harmonics of degree $n$. The corresponding Laplace eigenvalues are $n(n+d-1)$, and hence the first eigenspace $E_1=\mathscr H_1$ has dimension $d+1$. For a smooth real mean-zero perturbation $\eta$, we can then write its degree-two component as $\eta_2(x)=x^{\mathsf T}Q_\eta x$, where $Q_\eta$ is a real symmetric trace-free matrix of size $d+1$, which may vanish even when $\eta\neq 0$.

\begin{introsphere}
On the unit round $S^d$, $d\ge2$, the inequality
$\mathfrak h_\rho\ge d\,\mathfrak g_\rho$ holds for a smooth
positive probability density $\rho$ if and only if $\rho\equiv1$.
Moreover, for every nonzero smooth real mean-zero $\eta$ and every
sufficiently small nonzero $\varepsilon$ with
$1+\varepsilon\eta>0$,
\begin{equation}\label{eq:intro-sphere-index}
\begin{aligned}
\operatorname{ind}_-(\mathfrak d_{1+\varepsilon\eta})
&=n_-(\varepsilon Q_\eta)+\dim\ker Q_\eta,\\
\operatorname{nul}(\mathfrak d_{1+\varepsilon\eta})&=0.
\end{aligned}
\end{equation}
\end{introsphere}

The rigidity statement holds for every smooth positive density, whereas the index formula describes the local behavior along each fixed perturbation direction, including arbitrary mixtures of harmonic degrees and infinite smooth expansions. In particular, even when $\eta_2=0$ and the first-order defect vanishes on $E_1$, the coupling with higher modes produces the maximal index $d+1$. The proof of Theorem~\ref{thm:sphere-index} uses an integration-by-parts identity on hemispheres, and the Gram calculation on the two-sphere gives explicit second-order eigenvalue shifts.

Although spheres and complete flags provide two important examples, they do not cover all the possibilities described by the general index formula. On $S^1\times S^2$, for instance, a nonconstant density mode can be spectrally resonant and preserve the sharp inequality locally, showing that the rigidity found on spheres does not extend to every homogeneous space. This difference is explained by the representation channels that can interact, the resonances that determine which couplings remain, and the multiplicities that enter the index. 

\subsection*{Organization of the paper}

We begin in Section~\ref{sec:preliminaries} with the notation for Wasserstein geometry and representation theory, before establishing the homogeneous spectral splitting and the triple-product formulas in Sections~\ref{sec:static} and~\ref{sec:exact-couplings}. These results lead to the general index reduction, the Gram operator, and the finite-channel test developed in Sections~\ref{sec:obstruction} and~\ref{sec:coercive}. We then apply the theory to round spheres in Section~\ref{sec:spheres}, where we prove spherical rigidity and the local classification, and to complete flags in Section~\ref{sec:complete-flags}, where we compute the corresponding signatures. A Fourier check is given in Section~\ref{sec:circle-example}. Finally, Appendix~\ref{sec:su2-calculation} presents the normalized $SU(2)$ intertwiners and Haar triple integrals, which are used in Appendix~\ref{sec:su2-detection} to give a second description of the spherical results through representation theory.

\section{Wasserstein and representation-theoretic preliminaries}
\label{sec:preliminaries}
\subsection{Compact homogeneous spaces and Peter--Weyl coefficients}
\label{sec:groups-potentials}

Let $G$ be a compact Lie group and $H\subset G$ a closed subgroup. Assume that $F=G/H$ is connected and has positive dimension, and equip it with a $G$-invariant Riemannian metric.
Writing $\Psi:G\to F$ for the quotient map and $dg$ for the normalized
Haar measure of \(G\), the unique invariant probability measure on $F$ is
\begin{equation}\label{normalizedmeasureF}
m=\Psi_\#dg=\frac{d\operatorname{vol}_g}{\operatorname{vol}_g(F)}.
\end{equation}

We consider the usual Hermitian inner product of \(L^2(F,m)\) linear in the first argument.
Functions on $F$ identify with right $H$-invariant functions on
$G$, and by the Peter-Weyl theorem \cite{BrockerDieck1985,Farashahi2017} 
\begin{equation}\label{L2homogeneous}
L^2(F,m)=\widehat{\bigoplus}_{\pi\in\widehat G}\mathcal W_\pi,
\qquad \mathcal W_\pi\simeq V_\pi^*\otimes V_\pi^H.
\end{equation}
where $\pi:G\to U(V_\pi)$ is an irreducible unitary representation of \(G\), $d_\pi=\dim V_\pi$,
and $r_\pi=\dim V_\pi^H$ is its multiplicity.
For orthonormal bases $v_{\pi,i}$ of $V_\pi$ and
$w_{\pi,a}$ of $V_\pi^H$, define the matrix coefficients
\[
\Phi_{\pi,a,i}(gH)=\sqrt{d_\pi}\,
\langle\pi(g)w_{\pi,a},v_{\pi,i}\rangle.
\]

By Schur orthogonality relations these functions are orthonormal
\begin{equation}\label{SchurHomogeneous}
\int_F\Phi_{\pi,a,i}\overline{\Phi_{\sigma,b,j}}\,dm
 =\delta_{\pi\sigma}\delta_{ab}\delta_{ij}.
\end{equation}

In particular, these functions form an orthonormal basis of
$L^2(F,m)$, and $\mathcal W_\pi$ consists of $r_\pi$ equivalent
copies of the irreducible $G$-module $V_\pi^*$.

\subsection{Otto metric and the entropy Hessian}
\label{sec:wasserstein-definitions}

Since $F$ is compact, we have
$\mathcal P_2(F)=\mathcal P(F)$. In the smooth Otto calculus,
the velocity of a curve $\mu_t=\rho_t m$ is represented
by a gradient field $\nabla\varphi_t$, where the potential
$\varphi_t$ is related to the density through the
continuity equation
\begin{equation}\label{eq:continuity-reader}
\partial_t\rho_t
+\operatorname{div}(\rho_t\nabla\varphi_t)=0.
\end{equation}
This representation of velocities gives the tangent space
at $\mu=\rho m$ as the $L^2(\mu)$-closure of smooth
gradient fields,
\[
T_\mu\mathcal P(F)
=\overline{\{\nabla\varphi:\varphi\in C^\infty(F)\}}
^{L^2(F,\mu)},
\]
with the Otto metric
\[
\langle\nabla\varphi,\nabla\psi\rangle_{W,\mu}
=\int_F\rho\langle\nabla\varphi,\nabla\psi\rangle\,dm.
\]
For a smooth Wasserstein geodesic, the velocity potential
also satisfies the Hamilton--Jacobi equation
\begin{equation}\label{eq:hamilton-jacobi}
\partial_t\varphi_t+\frac12|\nabla\varphi_t|^2=0.
\end{equation}

We study the entropy functional associated with the
reference measure $m$,
\begin{equation}\label{eq:entropy}
\Ent_m(\rho m)=\int_F\rho\log\rho\,dm.
\end{equation}
Its second variation along Wasserstein geodesics
is described by the Wasserstein Hessian, which,
for smooth potentials $\varphi,\psi$, is given by
\begin{equation}\label{eq:hessian-general-density}
\Hess^W\Ent_m\big|_{\rho m}(\nabla\varphi,\nabla\psi)
=\int_F\rho\bigl(
\langle\Hess\varphi,\Hess\psi\rangle
+\Ric(\nabla\varphi,\nabla\psi)
\bigr)\,dm.
\end{equation}
These formulas follow the smooth Otto calculus
\cite{Otto2001,Lott2008} and its
Bochner--Bakry--\'Emery background
\cite{Bakry1991}, see also \cite{Villani2009}.

At the reference measure $m$, the integrated Bochner
identity gives a simpler expression for this Hessian,
which will be used throughout the paper.

\begin{lemma}\label{lem:bochner-hessian}
For smooth complex-valued potentials $\varphi,\psi$,
\begin{equation}\label{eq:integrated-bochner}
\Hess^W\Ent_m\big|_m(\nabla\varphi,\nabla\psi)
=\int_F(\Delta\varphi)\overline{\Delta\psi}\,dm.
\end{equation}
\end{lemma}

\begin{proof}
The result follows by integrating the polarized
Bochner identity over the closed manifold and
applying integration by parts.
\end{proof}

\subsection{Spectral realization of the entropy Hessian}
\label{sec:spectral-realization}

We use the following form construction to define the generalized Hessian eigenvalues, without assuming homogeneity. The spaces may also be complexified using the Hermitian convention introduced above, while indices and nullities are always taken over the reals.

\begin{proposition}\label{prop:spectral-realization}
Let $F$ be a closed connected Riemannian manifold of positive
dimension, let $m$ be the normalized volume, and let $\rho$ be a smooth
positive probability density. Set
\[
X=H^1(F)\cap L^2_0(F,m),\qquad
\mathcal V=H^2(F)\cap L^2_0(F,m).
\]
The form $\mathfrak g_\rho$ is an inner product on $X$ equivalent
to its $H^1$ norm. On this Hilbert space, $\mathfrak h_\rho$ is a
densely defined closed form with domain $\mathcal V$, bounded from
below. Its associated self-adjoint operator $A_\rho$ has compact
resolvent and smooth eigenfunctions, characterized by
\begin{equation}\label{eq:hessian-realization}
\mathfrak h_\rho(u,v)=\kappa\mathfrak g_\rho(u,v),
\qquad v\in\mathcal V.
\end{equation}
For each real $a$, the negative index of
$\mathfrak h_\rho-a\mathfrak g_\rho$ equals the number of
eigenvalues of $A_\rho$ below $a$, counted with multiplicity,
its nullity equals the multiplicity at $a$.

For $\rho_\varepsilon=1+\varepsilon\eta$, with $\eta$ smooth,
real, and mean-zero, the forms admit a real-analytic realization
on a fixed Hilbert space near $\varepsilon=0$, with common form
domain $\mathcal V$. Every isolated finite eigenvalue cluster
therefore has the analytic perturbation theory of self-adjoint
forms.
\end{proposition}
\begin{proof}
The bounds $0<\min\rho\leq\max\rho<\infty$ and the Poincar\'e
inequality give the asserted norm equivalence for $\mathfrak g_\rho$.
Boundedness of the curvature and the standard second-derivative
Sobolev estimate give constants $a_\rho>0$ and $b_\rho$ such that
\[
\mathfrak h_\rho(u,u)
\geq a_\rho\|u\|_{H^2}^2-b_\rho\|u\|_{H^1}^2,
\qquad u\in\mathcal V.
\]
Indeed, the principal term bounds
$(\min\rho)\|\Hess u\|_{L^2}^2$ from below and the Ricci term
is bounded by a constant times $\|\nabla u\|_{L^2}^2$.
For sufficiently large $c$, the form
$\mathfrak h_\rho+c\mathfrak g_\rho$ is thus coercive and
continuous on $\mathcal V$, with norm equivalent to $H^2$.
It follows that $\mathfrak h_\rho$ is closed and semibounded.
The representation theorem for closed forms
\cite[Ch.~VI]{Kato1995} defines $A_\rho$. Its shifted resolvent
maps $X$ boundedly into $\mathcal V$ by the coercive estimate.
The compact embedding $\mathcal V\hookrightarrow X$ makes this
resolvent compact. The spectral theorem and the min--max principle
give the assertions about index and nullity. Smooth potentials
have the same index, by density in $\mathcal V$.
The weak eigenvalue equation is a fourth-order elliptic equation
whose principal symbol is $\rho(x)|\xi|_g^4$, elliptic regularity
therefore makes its eigenfunctions smooth.

For the analytic assertion use $\mathfrak g_1$ as fixed inner
product on $X$. Define the bounded self-adjoint operator $T$ by
\[
\mathfrak g_1(Tu,v)=\int_F\eta\langle\nabla u,\nabla v\rangle\,dm.
\]
It also acts boundedly on $\mathcal V$: for $u\in\mathcal V$,
$LTu=-\operatorname{div}(\eta\nabla u)\in L^2_0$, and the inverse
of $L=-\Delta$ on mean-zero functions maps $L^2_0$ into $H^2$.
Let $G_\varepsilon=I+\varepsilon T$. For small $|\varepsilon|$,
its powers $G_\varepsilon^{\pm1/2}$ have norm-convergent binomial
series both on $X$ and on $\mathcal V$.
The map $U_\varepsilon=G_\varepsilon^{1/2}$ is unitary from
$(X,\mathfrak g_{\rho_\varepsilon})$ to $(X,\mathfrak g_1)$,
preserves $\mathcal V$, and transforms the Hessian into
\[
\widetilde{\mathfrak h}_\varepsilon(p,q)
=\mathfrak h_{\rho_\varepsilon}
 (G_\varepsilon^{-1/2}p,G_\varepsilon^{-1/2}q).
\]
This is a real-analytic family of closed symmetric forms with
common domain $\mathcal V$ and locally uniform coercive shifts.
Analytic form perturbation theory \cite[Ch.~VII]{Kato1995}
applies to its isolated eigenvalue clusters. Unitary equivalence
returns the same eigenvalues for \eqref{eq:hessian-realization}.
\end{proof}

\section{The entropy Hessian at the invariant measure}\label{sec:static}

The invariant Laplacian commutes with the left $G$-action, so its
spectral decomposition is controlled by the multiplicity spaces in
Peter--Weyl theory. More precisely, Schur's lemma reduces its
restriction to each isotypical summand to a self-adjoint operator
on $V_\pi^H$ \cite{BrockerDieck1985}. This explains why
multiplicities enter the Wasserstein calculation already at the
homogeneous measure.

\begin{lemma}
\label{lem:schur-laplacian}
For each $\pi\in\widehat G$, the Peter--Weyl component
$\mathcal W_\pi$ is invariant under $-\Delta$.
Moreover, under the unitary identification
$\mathcal W_\pi\simeq V_\pi^*\otimes V_\pi^H$,
there exists a unique nonnegative self-adjoint operator
$A_\pi:V_\pi^H\to V_\pi^H$ such that
\begin{equation}\label{eq:schur-block}
-\Delta\big|_{\mathcal W_\pi}
=
I_{V_\pi^*}\otimes A_\pi.
\end{equation}
\end{lemma}

\begin{proof}
The invariant metric makes $-\Delta$ commute with the left
$G$-action and preserve each isotypical summand.
The commutant on $V_\pi^*\otimes V_\pi^H$ is
$I_{V_\pi^*}\otimes\operatorname{End}(V_\pi^H)$ by Schur's
lemma. This gives $A_\pi$ and its uniqueness.
Self-adjointness and nonnegativity of $-\Delta$ transfer to
$A_\pi$ by testing on vectors $v\otimes w$ with $\|v\|=1$.
\end{proof}

Choose an orthonormal eigenbasis $A_\pi w_{\pi,a}=\lambda_{\pi,a}w_{\pi,a}$ and define
\begin{equation}\label{eq:spectral-subblock}
\mathcal W_{\pi,a}=V_\pi^*\otimes\operatorname{Span}\{w_{\pi,a}\},
\qquad
\mathcal E_{\pi,a}=\{\nabla\varphi:\ \varphi\in\mathcal W_{\pi,a}\}.
\end{equation}
For repeated eigenvalues, grouping equal eigendirections gives the
basis-independent version of this notation. The point of this decomposition is
that, at the homogeneous measure $m$, both the Otto metric and the entropy
Hessian are diagonalized by gradients of Laplace eigenfunctions, so that on each
subblock the Hessian is a scalar multiple of the metric.

\begin{proposition}
\label{thm:main-homogeneous}
The complexified Otto tangent space admits the orthogonal Hilbert sum
\begin{equation}\label{eq:tangent-decomposition}
T_m\mathcal P_2(F)_{\mathbb C}
=\widehat{\bigoplus_{\substack{\pi\in\widehat G,\ a\\ \lambda_{\pi,a}>0}}}
\mathcal E_{\pi,a}.
\end{equation}
If $V=\sum_{\pi,a}V_{\pi,a}$ is in the form domain of the entropy Hessian, then
\begin{equation}\label{eq:global-hessian-spectral}
\operatorname{Hess}^W\!\operatorname{Ent}_m(V,V)
=\sum_{\pi,a}\lambda_{\pi,a}\,\|V_{\pi,a}\|_{W,m}^2.
\end{equation}
Consequently, on every nonzero subblock,
\begin{equation}\label{eq:Rayleigh-block}
\frac{\operatorname{Hess}^W\!\operatorname{Ent}_m(V,V)}{\|V\|_{W,m}^2}
=\lambda_{\pi,a},
\qquad 0\neq V\in\mathcal E_{\pi,a}.
\end{equation}
Thus $\lambda_{\pi,a}$ is the infinitesimal Hessian rate at the homogeneous
measure, equivalently the generalized Hessian eigenvalue relative to the Otto
metric on that subblock.
\end{proposition}

\begin{proof}
Let $u\in\mathcal W_{\pi,a}$ and $v\in\mathcal W_{\sigma,b}$
be Laplace eigenfunctions with eigenvalues $\lambda_{\pi,a}$ and
$\lambda_{\sigma,b}$. Integration by parts and
Lemma~\ref{lem:bochner-hessian} give
\begin{equation}\label{eq:Otto-eigen-norm}
\mathfrak g_1(u,v)=\lambda_{\pi,a}\langle u,v\rangle_{L^2(m)},
\qquad
\mathfrak h_1(u,v)=\lambda_{\pi,a}\lambda_{\sigma,b}
\langle u,v\rangle_{L^2(m)}.
\end{equation}
The eigenspaces are mutually orthogonal. Since the only
zero-eigenvalue functions are constants, their gradients vanish.
Finite spectral sums are dense in the $H^1$ topology and hence
their gradients are dense in the Otto tangent space. This proves
\eqref{eq:tangent-decomposition} and Parseval's identity.
On each nonconstant summand the Hessian equals
$\lambda_{\pi,a}$ times the Otto metric, summing the diagonal
terms proves \eqref{eq:global-hessian-spectral} on finite sums.
Approximation in the form norm yields the asserted identity on
the full form domain, and restriction to one summand gives
\eqref{eq:Rayleigh-block}.
\end{proof}

Thus the spectrum of the Laplacian determines the Hessian at $m$ completely,
and its first positive eigenvalue yields a uniform lower bound.

\begin{corollary}\label{cor:static-gap}
If $\lambda_1(F)>0$ is the first positive eigenvalue of $-\Delta$, then for
every $V$ in the form domain,
\begin{equation}\label{eq:global-static-bound}
\operatorname{Hess}^W\!\operatorname{Ent}_m(V,V)\geq\lambda_1(F)\,\|V\|_{W,m}^2.
\end{equation}
The equality directions are precisely the gradients of first
Laplace eigenfunctions.
\end{corollary}

\begin{proof}
Every eigenvalue in \eqref{eq:global-hessian-spectral} is at least
$\lambda_1(F)$. Apply this bound termwise and use orthogonality
of the tangent-space decomposition. Equality is attained on $E_1$.
\end{proof}

\section{Coupling formulas}
\label{sec:exact-couplings}

Let
\begin{equation}\label{eq:rho-assumptions}
\mu=\rho m,\qquad \rho\in C^\infty(F,\mathbb R),\qquad
0<c\leq\rho\leq C,\qquad \int_F\rho\,dm=1.
\end{equation}
The bounds give equivalent weighted Sobolev norms. Define
\begin{align}
\mathfrak g_\rho(\varphi,\psi)
&=\int_F\rho\langle\nabla\varphi,\nabla\psi\rangle\,dm,
\label{eq:g-rho}\\
\mathfrak h_\rho(\varphi,\psi)
&=\int_F\rho\big(\langle\Hess\varphi,\Hess\psi\rangle
+\Ric(\nabla\varphi,\nabla\psi)\big)\,dm.
\label{eq:h-rho}
\end{align}
They are the complexified Otto metric and entropy Hessian at $\mu$.

Let $\{\Phi_A\}_{A\geq0}$ be an orthonormal eigenbasis adapted to the
subblocks of Proposition~\ref{thm:main-homogeneous}:
\begin{equation}\label{eq:eigenbasis}
-\Delta\Phi_A=\lambda_A\Phi_A,\qquad
\int_F\Phi_A\overline{\Phi_B}\,dm=\delta_{AB},\qquad \Phi_0=1.
\end{equation}
The multi-index $A=(\pi,a,i)$ records representation type,
multiplicity direction, and basis vector. Since $\rho$ is real and
smooth,
\begin{equation}\label{eq:density-expansion}
\rho=1+\eta,\qquad
\eta=\sum_{C\geq1}r_C\Phi_C,\qquad
r_C=\int_F\eta\overline{\Phi_C}\,dm,
\end{equation}
with convergence in every Sobolev norm and in $C^\infty$. Reality
imposes conjugate symmetry on the coefficients in a basis closed under
conjugation.

A Clebsch--Gordan coefficient is a matrix entry of an isometric
intertwiner $J:V_\sigma\to V_\pi\otimes V_\tau$, see
\cite{FultonHarris1991}. The formulas below combine these entries
with Laplace eigenvalues.

\begin{proposition}\label{prop:exact-coupling}
Assume \eqref{eq:rho-assumptions}--\eqref{eq:density-expansion}
and define the triple-product coefficients
\begin{equation}\label{eq:triple-C}
\mathcal C_{ABC}
=\int_F\Phi_A\overline{\Phi_B}\Phi_C\,dm.
\end{equation}
Then the metric and the entropy Hessian admit the exact
spectral expansions
\begin{align}
\mathfrak g_\rho(\Phi_A,\Phi_B)
&=\lambda_A\delta_{AB}
+\frac12\sum_{C\geq1}r_C
(\lambda_A+\lambda_B-\lambda_C)\mathcal C_{ABC},
\label{eq:exact-G}\\
\mathfrak h_\rho(\Phi_A,\Phi_B)
&=\lambda_A\lambda_B\delta_{AB}
+\frac14\sum_{C\geq1}r_C
(\lambda_A+\lambda_B-\lambda_C)^2\mathcal C_{ABC}.
\label{eq:exact-H}
\end{align}
Both series converge absolutely.
\end{proposition}

\begin{proof}
The product rule and the eigenvalue equations give
\begin{equation}\label{eq:product-rule-Hermitian}
\Delta(\Phi_A\overline{\Phi_B})
=-(\lambda_A+\lambda_B)\Phi_A\overline{\Phi_B}
+2\langle\nabla\Phi_A,\nabla\Phi_B\rangle.
\end{equation}
Multiplying by $\Phi_C$ and integrating by parts yields
\begin{equation}\label{eq:gradient-triple}
\int_F\Phi_C\langle\nabla\Phi_A,\nabla\Phi_B\rangle\,dm
=\frac{\lambda_A+\lambda_B-\lambda_C}{2}\mathcal C_{ABC}.
\end{equation}
Combining this with $\mathfrak g_1(\Phi_A,\Phi_B)
=\lambda_A\delta_{AB}$ proves \eqref{eq:exact-G}.

The polarized Bochner identity applied to eigenfunctions gives
\begin{equation}\label{eq:pointwise-Bochner-eigen}
\begin{aligned}
&\langle\Hess\Phi_A,\Hess\Phi_B\rangle
+\Ric(\nabla\Phi_A,\nabla\Phi_B)\\
&\quad=\frac12\Delta\langle\nabla\Phi_A,\nabla\Phi_B\rangle
+\frac{\lambda_A+\lambda_B}{2}
\langle\nabla\Phi_A,\nabla\Phi_B\rangle.
\end{aligned}
\end{equation}
Integrating against $\Phi_C$, and using
\eqref{eq:gradient-triple}, gives
\[
\int_F\Phi_C\bigl(\langle\Hess\Phi_A,\Hess\Phi_B\rangle
+\Ric(\nabla\Phi_A,\nabla\Phi_B)\bigr)\,dm
=\frac{(\lambda_A+\lambda_B-\lambda_C)^2}{4}
\mathcal C_{ABC}.
\]
The homogeneous Hessian entry equals
$\lambda_A\lambda_B\delta_{AB}$ by
\eqref{eq:integrated-bochner}, proving \eqref{eq:exact-H}.

For completeness, smoothness of $\eta$ implies, for every
positive integer $N$,
\[
|r_C|=|\langle\eta,\Phi_C\rangle|
\le(1+\lambda_C)^{-N}\|(I-\Delta)^N\eta\|_{L^2(m)}.
\]
For fixed $A,B$, Cauchy--Schwarz gives
$|\mathcal C_{ABC}|\le\|\Phi_A\overline{\Phi_B}\|_{L^2}$.
Thus the summands grow at most quadratically in $\lambda_C$
while $r_C$ decays faster than any inverse power.  Weyl's law
implies absolute convergence for $N>\frac12\dim F+2$,
Sobolev convergence justifies the integrations above.
\end{proof}

The coefficients $\mathcal C_{ABC}$ describe how different
spectral blocks interact. However, not every combination
of representation types can contribute. The following
selection rule identifies a necessary condition for
a nonzero coupling.

\begin{lemma}\label{lem:selection-rule}
Let $\Phi_A$, $\Phi_B$, and $\Phi_C$ have representation
types $\pi$, $\sigma$, and $\tau$, respectively.
If $\mathcal C_{ABC}\neq0$, then
\begin{equation}\label{eq:fusion-condition}
N_{\pi,\sigma^*}^{\tau^*}>0,
\qquad
N_{p,q}^{r}
:=\dim\Hom_G(V_r,V_p\otimes V_q).
\end{equation}
\end{lemma}

\begin{proof}
Since $m$ is $G$-invariant, the change of variables
$y=g^{-1}x$ gives
\[
\begin{aligned}
&\int_F (g\cdot\Phi_A)
\overline{(g\cdot\Phi_B)}
(g\cdot\Phi_C)\,dm\\
&\qquad=
\int_F\Phi_A(y)\overline{\Phi_B(y)}
\Phi_C(y)\,dm(y)
=\mathcal C_{ABC}.
\end{aligned}
\]
Thus, integration defines a $G$-invariant trilinear
functional on the corresponding representation spaces.
Equivalently, it determines an intertwining map
\[
V_\pi\otimes V_{\sigma^*}
\longrightarrow V_{\tau^*}.
\]
If $\mathcal C_{ABC}\neq0$, this map is nonzero. Now, since $G$ is a compact Lie group, the tensor product decomposes as
\[
V_\pi\otimes V_{\sigma^*}
\simeq
\bigoplus_{r\in\widehat G}
N_{\pi,\sigma^*}^{r}V_r.
\]
By Schur's lemma, a nonzero intertwining map into
$V_{\tau^*}$ exists only if $V_{\tau^*}$ occurs
in this decomposition. Therefore,
\[
N_{\pi,\sigma^*}^{\tau^*}>0.
\]
\end{proof}
The previous lemma identifies which representation types
may interact. We now use the Clebsch--Gordan decomposition
to compute the corresponding coupling coefficients explicitly,
taking into account the $H$-fixed multiplicities.
\begin{lemma}
\label{lem:CG-explicit}
Choose orthonormal eigenvectors $w_{\pi,a}\in V_\pi^H$ of
$A_\pi$ and normalized coefficient functions as in
\eqref{SchurHomogeneous}. Here
$\mathcal C_{ABC}=\int\Phi_A\overline{\Phi_B}\Phi_C\,dm$
contains a conjugate in its second factor, even when the
Peter--Weyl basis is complex. Equivalently, the selection condition
may be written as either
$\operatorname{Hom}_G(V_\sigma,V_\pi\otimes V_\tau)\ne0$ or
$\operatorname{Hom}_G(V_\pi\otimes V_\sigma^*,V_\tau^*)\ne0$.
For a triple of types
$\pi,\tau,\sigma$, let
\[
J_\nu:V_\sigma\longrightarrow V_\pi\otimes V_\tau,
\qquad \nu=1,\ldots,N_{\pi\tau}^{\sigma},
\]
be mutually orthogonal isometric $G$-intertwiners whose images
span the $\sigma$-isotypical summand of $V_\pi\otimes V_\tau$.
For $A=(\pi,a,i)$, $C=(\tau,c,k)$, and $B=(\sigma,b,j)$,
\begin{equation}\label{eq:CG-H-fixed}
\begin{aligned}
\mathcal C_{ABC}
={}&\sqrt{\frac{d_\pi d_\tau}{d_\sigma}}
\sum_{\nu=1}^{N_{\pi\tau}^{\sigma}}
\langle w_{\pi,a}\otimes w_{\tau,c},
       J_\nu w_{\sigma,b}\rangle\cdot
\langle J_\nu v_{\sigma,j},
      v_{\pi,i}\otimes v_{\tau,k}\rangle.
\end{aligned}
\end{equation}
The Hermitian products are linear in their first variable, as fixed
in \eqref{SchurHomogeneous}, when $N_{\pi\tau}^{\sigma}=0$ the sum
is zero.  In particular, tensor-product multiplicity alone does not
settle nonvanishing: the $H$-fixed-vector contractions in
\eqref{eq:CG-H-fixed} must also be evaluated.
\end{lemma}

\begin{proof}
Set $U=\pi\otimes\tau$,
$w=w_{\pi,a}\otimes w_{\tau,c}$, and
$v=v_{\pi,i}\otimes v_{\tau,k}$. Multiplication of matrix
coefficients gives
\[
\mathcal C_{ABC}
=\sqrt{d_\pi d_\tau d_\sigma}
\int_G\langle U(g)w,v\rangle
\overline{\langle\sigma(g)w_{\sigma,b},v_{\sigma,j}\rangle}\,dg.
\]
Only the $\sigma$-isotypical component of $U$ contributes, by
Schur orthogonality. Its orthogonal projector is
$P_\sigma=\sum_\nu J_\nu J_\nu^*$, and equivariance yields
\[
\langle U(g)P_\sigma w,v\rangle
=\sum_\nu\langle\sigma(g)J_\nu^*w,J_\nu^*v\rangle.
\]
With the Hermitian pairing linear in its first argument,
\[
\int_G\langle\sigma(g)a,b\rangle
\overline{\langle\sigma(g)c,d\rangle}\,dg
=\frac{\langle a,c\rangle\langle d,b\rangle}{d_\sigma}.
\]
Applying this identity termwise gives
\begin{equation}\label{eq:CG-four-steps}
\mathcal C_{ABC}
=\sqrt{\frac{d_\pi d_\tau}{d_\sigma}}
\sum_\nu
\langle w,J_\nu w_{\sigma,b}\rangle
\langle J_\nu v_{\sigma,j},v\rangle,
\end{equation}
which is \eqref{eq:CG-H-fixed}.
\end{proof}

For $H=\{e\}$, choose orthonormal bases and set
$C^\nu_{ik,j}=\langle J_\nu e_j^\sigma,
e_i^\pi\otimes e_k^\tau\rangle$. The coordinate form is
\begin{equation}\label{eq:CG-coordinate-contraction}
\int_G\Phi_{\pi,a,i}\Phi_{\tau,c,k}
\overline{\Phi_{\sigma,b,j}}\,dg
=\sqrt{\frac{d_\pi d_\tau}{d_\sigma}}
\sum_\nu\overline{C^\nu_{ac,b}}\,C^\nu_{ik,j}.
\end{equation}
For general $H$, the first contraction uses the $H$-fixed
multiplicity vectors in \eqref{eq:CG-H-fixed}. Thus the permitted
representation type is only a necessary condition for an active
channel: the contractions must be nonzero after summation.

\section{Spectral channels and the index of the entropy-Hessian defect}
\label{sec:obstruction}

The coupling formulas now allow us to count the directions in
which the entropy Hessian falls below its sharp homogeneous bound.
Throughout this section, $L=-\Delta$, $\lambda_1$ is its first
positive eigenvalue, and $\mathscr E_\alpha=\ker(L-\alpha I)$, with
$\mathscr E_\alpha=\{0\}$ if $\alpha$ is not in the spectrum.
We reserve $E_1=\mathscr E_{\lambda_1}$ for the first positive
eigenspace, the subscript of $\mathscr E_\alpha$ is always a
numerical eigenvalue.
\begin{definition}
\label{def:obstruction}
For a smooth positive probability density $\rho$, we define the entropy-Hessian defect by the bilinear form
\[
\mathfrak d_{\rho}(u,v)
=\mathfrak h_{\rho}(u,v)-\lambda_{1}\mathfrak g_{\rho}(u,v).
\]
This defect measures the failure of the entropy Hessian to satisfy the spectral lower bound $\mathfrak h_{\rho}\geq\lambda_1\mathfrak g_{\rho}$.

To quantify this failure, define $\operatorname{ind}_{-}(\mathfrak d_{\rho})$
as the maximal dimension of a real subspace of smooth potentials
modulo constants on which $\mathfrak d_{\rho}$ is negative definite.
The same index is obtained on $H^2(F,\mathbb R)/\mathbb R$:
the form is continuous in the $H^2$ topology, smooth potentials
are dense in $H^2$, and approximating a basis of any finite-dimensional
negative subspace preserves the negative definiteness of its Gram
matrix.  Its nullity is the dimension of the kernel of the
corresponding continuous form on $H^2(F)/\mathbb R$.
Equivalently, Proposition~\ref{prop:spectral-realization}
identifies $\operatorname{ind}_{-}(\mathfrak d_{\rho})$ with the
number of generalized eigenvalues $\kappa<\lambda_1$ of
$\mathfrak h_\rho-\kappa\mathfrak g_\rho$, counted with
multiplicity.
\end{definition}

At $\rho=1$, the integrated Bochner identity gives
$\mathfrak h_{1}(u,v)=\langle Lu,Lv\rangle_{L^{2}(m)}$ and
$\mathfrak g_{1}(u,v)=\langle Lu,v\rangle_{L^{2}(m)}$. Therefore,
\begin{equation}\label{eq:D-homogeneous}
\mathfrak d_{1}(u,v)
=\langle Lu,Lv\rangle_{L^{2}(m)}
 -\lambda_{1}\langle Lu,v\rangle_{L^{2}(m)},
\qquad u,v\in H^2(F)/\mathbb R.
\end{equation}
For a Laplace eigenfunction $Lu=\beta u$, it follows that
\[
\mathfrak d_{1}(u,u)=\beta(\beta-\lambda_{1})\|u\|_{L^{2}}^{2}.
\]
Thus, the homogeneous defect is always nonnegative for mean-zero potentials. It is zero on $E_1$, and the terms involving one potential in $E_1$ and another in $E_1^\perp$ also vanish.

Expanding the density as in \eqref{eq:density-expansion} and
subtracting $\lambda_1$ times \eqref{eq:exact-G} from
\eqref{eq:exact-H}, we obtain the spectral factor
\[
\begin{aligned}
&\frac{1}{4}(\lambda_{A}+\lambda_{B}-\lambda_{C})^{2}
 -\frac{\lambda_{1}}{2}(\lambda_{A}+\lambda_{B}-\lambda_{C})\\
&\qquad=\frac{1}{4}(\lambda_{A}+\lambda_{B}-\lambda_{C})
 (\lambda_{A}+\lambda_{B}-\lambda_{C}-2\lambda_{1}).
\end{aligned}
\]
Consequently,
\begin{equation}\label{eq:D-triple}
\begin{aligned}
\mathfrak d_{\rho}(\Phi_{A},\Phi_{B})
={}&\lambda_{A}(\lambda_{A}-\lambda_{1})\delta_{AB}\\
&+\frac{1}{4}\sum_{C\geq 1}r_{C}
 (\lambda_{A}+\lambda_{B}-\lambda_{C})
 (\lambda_{A}+\lambda_{B}-\lambda_{C}-2\lambda_{1})
 \mathcal C_{ABC}.
\end{aligned}
\end{equation}
A term can therefore vanish either because its representation
channel is absent or because the numerical eigenvalues satisfy a
resonance relation. The next lemma makes this distinction precise
for couplings involving the first eigenspace.

\begin{lemma}\label{lem:first-row}
Let $\gamma>0$ be a Laplace eigenvalue and let
$\eta\in \mathscr E_{\gamma}$ be real and mean-zero. Suppose that
$\rho_{\varepsilon}=1+\varepsilon\eta>0$.
If $u\in E_{1}$ and $v\in \mathscr E_{\beta}$, where
$\beta\geq0$ is a Laplace eigenvalue, then
\begin{equation}\label{eq:resonance-coupling}
\mathfrak d_{\rho_{\varepsilon}}(u,v)
=\frac{\varepsilon}{4}
\bigl[(\beta-\gamma)^{2}-\lambda_{1}^{2}\bigr]
\int_{F}\eta u\overline{v}\,dm.
\end{equation}
For $\varepsilon\neq0$, all first-row pairings vanish,
for every $u\in E_{1}$ and every Laplace eigenspace,
if and only if
\begin{equation}\label{eq:resonance-invisible}
\eta E_{1}\subset
\mathscr E_{\gamma+\lambda_{1}}\oplus \mathscr E_{\gamma-\lambda_{1}}.
\end{equation}
A summand on the right is understood to be $\{0\}$
if its parameter is negative or is not in the spectrum.
\end{lemma}

\begin{proof}
For $Lu=\lambda_1u$, $Lv=\beta v$, and $L\eta=\gamma\eta$,
the homogeneous pairing vanishes and \eqref{eq:D-triple} gives
\[
\mathfrak d_{1+\varepsilon\eta}(u,v)
=\frac{\varepsilon}{4}
(\lambda_1+\beta-\gamma)
(\beta-\gamma-\lambda_1)
\langle\eta u,v\rangle_{L^2(m)},
\]
which is \eqref{eq:resonance-coupling}.  Fix $u\in E_1$.
For every $v\in \mathscr E_\beta$, orthogonal spectral projection yields
\[
\langle\eta u,v\rangle
=\langle\Pi_\beta(\eta u),v\rangle.
\]
The pairing vanishes for all $v\in \mathscr E_\beta$ exactly when
$[(\beta-\gamma)^2-\lambda_1^2]\Pi_\beta(\eta u)=0$.
The scalar factor vanishes precisely for
$\beta=\gamma+\lambda_1$ or $\beta=\gamma-\lambda_1$.
Applying this observation to every eigenspace and every
$u\in E_1$ proves \eqref{eq:resonance-invisible} in both
directions.
\end{proof}

These two conditions play different roles in determining whether a coupling is nonzero. For a density mode $\eta_{\gamma}\in \mathscr E_{\gamma}$ and a target eigenspace $\mathscr E_{\beta}$, consider the projection map and the corresponding spectral factor
\[
\begin{aligned}
T_{\beta,\gamma}&:E_{1}\longrightarrow \mathscr E_{\beta},
\qquad T_{\beta,\gamma}(u)=\Pi_{\beta}(\eta_{\gamma}u),\\[4pt]
c_{\beta,\gamma}&=\frac{1}{4}
\bigl[(\beta-\gamma)^{2}-\lambda_{1}^{2}\bigr].
\end{aligned}
\]

The projection map determines whether the product of spectral
modes has a component in $\mathscr E_{\beta}$, while the
spectral factor determines whether this component contributes
to the defect. Thus, a nonzero coupling requires both to be
nonzero.

For $u,v\in E_{1}$, the pointwise Bochner formula is given by
\begin{equation}\label{eq:first-space-density-Laplacian}
\mathfrak d_\rho(u,v)
=\frac12\int_F (\Delta\rho)
\langle\nabla u,\nabla v\rangle\,dm.
\end{equation}
Indeed, the pointwise defect integrand on two first eigenfunctions
is $\tfrac12\Delta\langle\nabla u,\nabla v\rangle$,
and integration by parts, using the self-adjointness of
$\Delta$, gives the desired identity.

\begin{lemma}\label{lem:tracefree}
On a compact connected homogeneous manifold, the restriction of
$\mathfrak d_{\rho}$ to the real first eigenspace has zero trace
for every smooth probability density $\rho$. If this restriction
is not zero, it has both positive and negative directions.
\end{lemma}
\begin{proof}
Let $u_1,\ldots,u_d$ be a real $L^2(m)$-orthonormal basis of
$E_1$.  Transitivity and orthogonality of the $G$-action imply
that $\sum_j u_j(x)^2$ is constant, integrating shows it equals
$d$.  Applying $\Delta$ and using $\Delta u_j=-\lambda_1u_j$
gives $\sum_j|\nabla u_j|^2=\lambda_1d$ pointwise.
Consequently, by \eqref{eq:first-space-density-Laplacian},
\[
\sum_{j=1}^d\mathfrak d_\rho(u_j,u_j)
=\frac12\int_F(\Delta\rho)
  \sum_{j=1}^d|\nabla u_j|^2\,dm
=\frac{\lambda_1d}{2}\int_F\Delta\rho\,dm=0.
\]
A nonzero real symmetric trace-free form has eigenvalues of
both signs.
\end{proof}

\section{The coercive complement and the effective finite-dimensional form}\label{sec:coercive}

Until Corollary~\ref{cor:detect-criterion}, $F$ may be any closed
connected Riemannian manifold of positive dimension, and $m$
denotes normalized volume. The notation $E_1$ and
$\mathscr E_\alpha$ retains the spectral meaning fixed above.

Fix a smooth real-valued function $\eta$ with zero mean and consider the density perturbation $\rho_\varepsilon=1+\varepsilon\eta$. We denote by $\mathfrak d_\eta$ the first variation of the defect at the uniform density:
\[
\mathfrak d_\eta
=\left.\frac{d}{d\varepsilon}\right|_{\varepsilon=0}
\mathfrak d_{1+\varepsilon\eta}.
\]
Thus, $\mathfrak d_\eta$ is a variation form, rather than the defect associated with a probability density. Since the defect depends linearly on $\rho$, this variation gives the exact identity
\[
\mathfrak d_{1+\varepsilon\eta}
=\mathfrak d_1+\varepsilon\mathfrak d_\eta.
\]

The natural form domain is the real Sobolev space of mean-zero
potentials $\mathcal V=H^2(F,\mathbb R)\cap L_0^2(F,m)$.
Orthogonal projection onto $E_1$ splits this space as
\begin{equation}\label{eq:critical-decomposition}
\mathcal V=E_1\oplus\mathcal V_+,
\qquad \mathcal V_+=\mathcal V\cap E_1^\perp,
\end{equation}
where orthogonality is taken in $L^2(F,m)$. The first summand
is precisely the kernel of the unperturbed defect, while the
second contains only higher Laplace modes. On this complement
the unperturbed form is coercive, as the spectral gap shows.
Let $\lambda_*>\lambda_1$ denote the next distinct positive eigenvalue of the Laplace operator. For $v\in\mathcal V_+$, its spectral expansion contains only eigenvalues $\beta\geq\lambda_*$:
\[
v=\sum_{\beta\geq\lambda_*}\Pi_\beta v.
\]
Using this expansion in \eqref{eq:D-homogeneous}, we obtain
\begin{equation}\label{eq:coercive-complement}
\begin{aligned}
\mathfrak d_{1}(v,v)
&=\sum_{\beta\geq\lambda_{*}}
  \beta(\beta-\lambda_{1})\|\Pi_{\beta}v\|_{L^{2}}^{2}\\
&=\sum_{\beta\geq\lambda_{*}}
  \beta^{2}\left(1-\frac{\lambda_{1}}{\beta}\right)
  \|\Pi_{\beta}v\|_{L^{2}}^{2}\\
&\geq\left(1-\frac{\lambda_{1}}{\lambda_{*}}\right)
  \sum_{\beta\geq\lambda_{*}}
  \beta^{2}\|\Pi_{\beta}v\|_{L^{2}}^{2}\\
&=\left(1-\frac{\lambda_{1}}{\lambda_{*}}\right)
  \|Lv\|_{L^{2}}^{2}.
\end{aligned}
\end{equation}
Set $c_*=1-\lambda_1/\lambda_*>0$. By elliptic regularity and the spectral gap, there exists a constant $C_F>0$ such that
\[
\|v\|_{H^2}\leq C_F\|Lv\|_{L^2},
\qquad v\in\mathcal V_+.
\]
Combining this estimate with the lower bound for $\mathfrak d_1$, we obtain
\[
\mathfrak d_1(v,v)
\geq c_*\|Lv\|_{L^2}^2
\geq c_*C_F^{-2}\|v\|_{H^2}^2.
\]

To control the perturbed defect, we also need an estimate for $\mathfrak d_\eta$. From the definitions of the entropy Hessian and the Otto metric, we have
\[
\mathfrak d_{\eta}(v,w)=\int_F\eta
\bigl(
\langle\Hess v,\Hess w\rangle
+\Ric(\nabla v,\nabla w)
-\lambda_1\langle\nabla v,\nabla w\rangle
\bigr)\,dm.
\]

We estimate each term using the Cauchy--Schwarz inequality. For the Hessian term, we obtain
\[
\begin{aligned}
\left|\int_F\eta\langle\Hess v,\Hess w\rangle\,dm\right|
&\leq\|\eta\|_{L^\infty}
\int_F|\Hess v|\,|\Hess w|\,dm\\[6pt]
&\leq\|\eta\|_{L^\infty}
\|\Hess v\|_{L^2}\|\Hess w\|_{L^2}.
\end{aligned}
\]

Set $R_{\Ric}=\sup_{x\in F}|\Ric_x|$. The remaining terms
are bounded in absolute value by
\[
(R_{\Ric}+\lambda_1)\|\eta\|_{L^\infty}
\|\nabla v\|_{L^2}\|\nabla w\|_{L^2}.
\]

It follows that
\begin{equation}\label{eq:eta-form-bound}
|\mathfrak d_{\eta}(v,w)|
\leq K_F\|\eta\|_{L^{\infty}}
\|v\|_{H^{2}}\|w\|_{H^{2}},
\qquad v,w\in\mathcal V_{+},
\end{equation}
where $K_F$ depends only on the fixed metric and $\lambda_{1}$.

Let $C_{\varepsilon}:\mathcal V_{+}\to\mathcal V_{+}^{*}$
be the operator induced by the restriction of
$\mathfrak d_{1+\varepsilon\eta}$ to $\mathcal V_{+}$.
The preceding bounds show that
\[
\begin{aligned}
\langle C_{\varepsilon}v,v\rangle
&=\mathfrak d_{1}(v,v)+\varepsilon\mathfrak d_{\eta}(v,v)\geq\left(\frac{c_{*}}{C_{F}^{2}}
       -|\varepsilon|K_F\|\eta\|_{L^{\infty}}\right)
       \|v\|_{H^{2}}^{2}.
\end{aligned}
\]
Choose $\varepsilon_{0}>0$ so that
$\varepsilon_{0}K_F\|\eta\|_{L^{\infty}}
<c_{*}/(2C_{F}^{2}) $
(for $\eta=0$ no restriction is needed). For
$|\varepsilon|<\varepsilon_0$, it follows that
\[
\langle C_{\varepsilon}v,v\rangle
\geq c_{*}(2C_{F}^{2})^{-1}\|v\|_{H^{2}}^{2}.
\]

Thus, the operator $C_\varepsilon$ is bounded and coercive on $\mathcal V_+$. By the Lax--Milgram theorem, it admits a bounded inverse
\[
C_\varepsilon^{-1}:\mathcal V_+^*\longrightarrow\mathcal V_+.
\]
This inverse acts between the Sobolev form domain and its dual, rather than on $L^2$.

With the defect controlled on $\mathcal V_+$, we now
examine how it couples this space with the first
eigenspace $E_1$. The following lemma uses this decomposition to reduce the negative index and nullity of the full defect to an operator on the finite-dimensional space $E_1$.

\begin{definition}
\label{def:MB-operators}
We define the operators
$M_\eta:E_1\to E_1$ and
$B_\eta:E_1\to\mathcal V_+^*$ by
\begin{equation}\label{eq:MB-operators}
\begin{aligned}
\langle M_\eta u,w\rangle
&=\mathfrak d_\eta(u,w),
&&u,w\in E_1,\\[6pt]
\langle B_\eta u,v\rangle
&=\mathfrak d_\eta(u,v),
&&u\in E_1,\quad v\in\mathcal V_+.
\end{aligned}
\end{equation}
The real symmetric operator $M_\eta$ describes
the first-order coupling within $E_1$, while
$B_\eta$ describes the coupling between $E_1$
and the higher eigenspaces.
\end{definition}

Together with $C_\varepsilon$, which controls
the defect on $\mathcal V_+$, these operators
describe the defect relative to the decomposition
$\mathcal V=E_1\oplus\mathcal V_+$.
Constant potentials do not contribute to these forms.

\begin{lemma}\label{lem:schur}
For $|\varepsilon|<\varepsilon_{0}$, let
\begin{equation}\label{eq:effective-schur}
S_{\varepsilon}
=\varepsilon M_{\eta}
 -\varepsilon^{2}B_{\eta}^{*}C_{\varepsilon}^{-1}B_{\eta}
\quad\text{on }E_{1}.
\end{equation}
For $u\in E_{1}$ and $v\in\mathcal V_{+}$, one has
\begin{equation}\label{eq:square-identity}
\begin{aligned}
\mathfrak d_{1+\varepsilon\eta}(u+v,u+v)
={}&\left\langle C_{\varepsilon}
\left(v+\varepsilon C_{\varepsilon}^{-1}B_{\eta}u\right),
 v+\varepsilon C_{\varepsilon}^{-1}B_{\eta}u\right\rangle +\langle S_{\varepsilon}u,u\rangle.
\end{aligned}
\end{equation}
In particular,
$\operatorname{ind}_{-}(\mathfrak d_{1+\varepsilon\eta})
=n_{-}(S_{\varepsilon})$ and
$\operatorname{nul}(\mathfrak d_{1+\varepsilon\eta})
=\dim\ker S_{\varepsilon}$.
\end{lemma}
\begin{proof}
Since $\mathfrak d_1$ vanishes whenever one argument lies in
$E_1$, its block form on $E_1\oplus\mathcal V_+$ is
\[
\begin{pmatrix}
\varepsilon M_\eta&\varepsilon B_\eta^*\\
\varepsilon B_\eta&C_\varepsilon
\end{pmatrix}.
\]
Thus, for $u\in E_1$ and $v\in\mathcal V_+$,
\[
\mathfrak d_{1+\varepsilon\eta}(u+v,u+v)
=\varepsilon\langle M_\eta u,u\rangle
+2\varepsilon\langle B_\eta u,v\rangle
+\langle C_\varepsilon v,v\rangle.
\]
Set $w=v+\varepsilon C_\varepsilon^{-1}B_\eta u$.
Expansion of $\langle C_\varepsilon w,w\rangle$ gives
\[
\begin{aligned}
\langle C_\varepsilon w,w\rangle
={}&\langle C_\varepsilon v,v\rangle
+2\varepsilon\langle B_\eta u,v\rangle\\
&+\varepsilon^2\langle B_\eta^*C_\varepsilon^{-1}B_\eta u,u\rangle.
\end{aligned}
\]
Subtracting the last term proves \eqref{eq:square-identity}.
The change of coordinates $(u,v)\mapsto(u,w)$ is invertible,
and $C_\varepsilon$ is positive. Hence all negative and null
directions come from the finite-dimensional form $S_\varepsilon$.
\end{proof}

The classical Schur-complement reduction now yields a local
index formula. The proof uses the Laplace spectral gap and the
Bochner identity, so the theorem is stated for an arbitrary closed
connected manifold. The spectral identities of
Lemma~\ref{lem:first-row} remain valid in any Laplace eigenbasis:
their derivation uses only the product rule and Bochner's formula.
To state the spectral form, we fix the notation
for the components of the mixed operator. Let $\Pi_\alpha$
denote the $L^2(m)$-orthogonal projection onto the Laplace
eigenspace with numerical eigenvalue $\alpha$, and write
$\eta=\sum_{\gamma>0}\eta_\gamma$, where
$\eta_\gamma=\Pi_\gamma\eta$. For each numerical eigenvalue
$\beta>\lambda_1$, define the map $B_{\eta,\beta}:E_1\to \mathscr E_\beta$ by
\begin{equation}\label{eq:matrix-Bbeta}
B_{\eta,\beta}u
=\frac14\sum_{\gamma>0}
\bigl[(\beta-\gamma)^2-\lambda_1^2\bigr]
\Pi_\beta(\eta_\gamma u).
\end{equation}
The notation $\beta$ and $\gamma$ always refers to numerical
Laplace eigenvalues, not to representation labels. The operators
$M_\eta$ and $B_\eta$ were defined in
\eqref{eq:MB-operators}, while $C_0$ is the positive homogeneous
defect form on $\mathcal V_+$. Write $d=\dim E_1$.

\begin{theorem}\label{thm:GH-index}
Let $F$ be a closed connected Riemannian manifold of positive dimension, equipped with its normalized Riemannian volume, and let $\eta$ be a smooth real mean-zero perturbation. Then there exists $\varepsilon_0>0$, depending on the metric and on $\eta$, such that, for every $0<|\varepsilon|<\varepsilon_0$ with $\rho_\varepsilon>0$, we have
\begin{equation}\label{eq:intro-index}
\begin{aligned}
\operatorname{ind}_-(\mathfrak d_{\rho_\varepsilon})
&=n_-(\varepsilon M_\eta)
+\operatorname{rank}(B_\eta|_{\ker M_\eta}),\\
\operatorname{nul}(\mathfrak d_{\rho_\varepsilon})
&=\dim(\ker M_\eta\cap\ker B_\eta).
\end{aligned}
\end{equation}
In particular, the negative index is at most $\dim E_1$.
\end{theorem}

The operators in this reduction admit the following spectral
formulas, which will make the index computable in the applications.

\begin{proposition}\label{prop:spectral-gram}
Under the hypotheses of Theorem~\ref{thm:GH-index},
the spectral components of the mixed operator satisfy the
positive semidefinite Gram identity
\begin{equation}\label{eq:gram-PW}
\displaystyle
B_{\eta}^{*}C_{0}^{-1}B_{\eta}
=\sum_{\beta>\lambda_{1}}
 \frac{B_{\eta,\beta}^{*}B_{\eta,\beta}}
 {\beta(\beta-\lambda_{1})}.
\end{equation}
The series converges in operator norm on $E_{1}$ and the kernel
of the Gram operator is $\ker B_{\eta}$. Finally,
for $u,w\in E_{1}$,
\begin{equation}\label{eq:matrix-M}
\langle M_{\eta}u,w\rangle
=\frac{1}{4}\sum_{\gamma>0}
 \gamma(\gamma-2\lambda_{1})
 \langle\eta_{\gamma}u,w\rangle_{L^{2}(m)}.
\end{equation}
On a homogeneous space $F=G/H$, the entries in
\eqref{eq:matrix-Bbeta} and \eqref{eq:matrix-M} are weighted sums
of the Peter--Weyl triple coefficients $\mathcal C_{ABC}$ from
Proposition~\ref{prop:exact-coupling}.
\end{proposition}

\begin{proof}
We first identify the blocks of the variation.  For $u\in E_1$
and $v\in \mathscr E_\beta$, Lemma~\ref{lem:first-row} and the
orthogonal projections of $\eta_\gamma u$ give
\begin{equation}\label{eq:proof-spectral-first-row}
\begin{aligned}
\mathfrak d_\eta(u,v)
&=\frac14\sum_{\gamma>0}
\bigl[(\beta-\gamma)^2-\lambda_1^2\bigr]
\langle\eta_\gamma u,v\rangle_{L^2(m)}\\
&=\left\langle
\frac14\sum_{\gamma>0}
\bigl[(\beta-\gamma)^2-\lambda_1^2\bigr]
\Pi_\beta(\eta_\gamma u),v\right\rangle_{L^2(m)}.
\end{aligned}
\end{equation}
This identifies $B_{\eta,\beta}$ as in
\eqref{eq:matrix-Bbeta}.  When $\beta=\lambda_1$, the
spectral factor becomes
$(\lambda_1-\gamma)^2-\lambda_1^2
=\gamma(\gamma-2\lambda_1)$, giving
\eqref{eq:matrix-M}.  For fixed smooth $u$ and $\eta$,
integration by parts in the definition of $\mathfrak d_\eta$
transfers all derivatives of $v$ to smooth coefficients. Hence
$\mathfrak d_\eta(u,v)=\langle f_{\eta,u},v\rangle_{L^2(m)}$
with $f_{\eta,u}$ smooth.  The higher-mode representative of
$B_\eta u$ is therefore
\[
B_\eta u=\sum_{\beta>\lambda_1}B_{\eta,\beta}u,
\]
with convergence in $L^2(F,m)$.

On $\mathscr E_\beta$, $\beta>\lambda_1$, the unperturbed complement
acts as
$C_0=\beta(\beta-\lambda_1)I$.  Orthogonality of distinct
eigenspaces consequently yields, for $u,w\in E_1$,
\begin{equation}\label{eq:proof-gram-pairing}
\begin{aligned}
\langle B_\eta^*C_0^{-1}B_\eta u,w\rangle
&=\langle C_0^{-1}B_\eta u,B_\eta w\rangle\\
&=\sum_{\beta>\lambda_1}
\frac{\langle B_{\eta,\beta}u,B_{\eta,\beta}w\rangle_{L^2(m)}}
{\beta(\beta-\lambda_1)}.
\end{aligned}
\end{equation}
This proves \eqref{eq:gram-PW}.  Its quadratic form is
\[
\langle \mathcal G_0u,u\rangle
=\sum_{\beta>\lambda_1}
\frac{\|B_{\eta,\beta}u\|_{L^2(m)}^2}
{\beta(\beta-\lambda_1)},
\qquad \mathcal G_0=B_\eta^*C_0^{-1}B_\eta.
\]
All denominators are positive, Parseval's identity therefore gives
$\mathcal G_0\ge0$ and $\ker \mathcal G_0=\ker B_\eta$.
For clarity, the series also converges in operator norm.  Let
$\lambda_+>\lambda_1$ be the next distinct Laplace eigenvalue
and $(e_j)_{j=1}^d$ an orthonormal basis of $E_1$. Then
\begin{equation}\label{eq:proof-gram-summability}
\begin{aligned}
&\sum_{\beta>\lambda_1}
\frac{\|B_{\eta,\beta}^*B_{\eta,\beta}\|_{\mathrm{op}}}
{\beta(\beta-\lambda_1)}\\
&\quad\le\frac1{\lambda_+(\lambda_+-\lambda_1)}
\sum_{j=1}^d\sum_{\beta>\lambda_1}
\|B_{\eta,\beta}e_j\|_{L^2(m)}^2\\
&\quad=\frac{\sum_{j=1}^d\|B_\eta e_j\|_{L^2(m)}^2}
{\lambda_+(\lambda_+-\lambda_1)}<\infty.
\end{aligned}
\end{equation}
Here we used
$\|B_{\eta,\beta}^*B_{\eta,\beta}\|_{\mathrm{op}}
\le\|B_{\eta,\beta}\|_{\mathrm{HS}}^2$.

This proves the spectral formulas and the convergence and kernel
assertions. On a homogeneous space, their Peter--Weyl expressions
follow from Proposition~\ref{prop:exact-coupling}.
\end{proof}

\begin{proof}[Proof of Theorem~\ref{thm:GH-index}]
We compute the inertia of the effective form
$S_\varepsilon$ in \eqref{eq:effective-schur}. Put
$\mathcal G_\varepsilon=B_\eta^*C_\varepsilon^{-1}B_\eta$, where
$C_\varepsilon=C_0+\varepsilon C_\eta$.
The complementary coercivity estimate makes
$C_\varepsilon^{-1}$ uniformly bounded. Its resolvent identity is
\[
C_\varepsilon^{-1}-C_0^{-1}
=-\varepsilon C_\varepsilon^{-1}C_\eta C_0^{-1},
\]
so, by boundedness of $C_\eta$ and $B_\eta$,
\begin{equation}\label{eq:proof-gram-variation}
\|\mathcal G_\varepsilon-\mathcal G_0\|_{\mathrm{op}}
\le C|\varepsilon|.
\end{equation}

Write
\[
\mathcal N=\ker M_\eta,\qquad R=\mathcal N^\perp,\qquad
\mathcal N_0=\mathcal N\cap\ker B_\eta,\qquad \mathcal N_1=\mathcal N\cap \mathcal N_0^\perp.
\]
Then $E_1=R\oplus \mathcal N_1\oplus \mathcal N_0$ orthogonally,
$M_R=M_\eta|_R$ is invertible, and
$\mathcal G_\varepsilon$ vanishes on $\mathcal N_0$ for every sufficiently small
$\varepsilon$.  Therefore $S_\varepsilon$ vanishes identically
on $\mathcal N_0$.  On the other hand, $\mathcal G_0$ is strictly positive on
$\mathcal N_1$: if $z\in \mathcal N_1\setminus\{0\}$, then $B_\eta z\ne0$,
and \eqref{eq:proof-gram-pairing} gives
$\langle \mathcal G_0z,z\rangle>0$.
Finite dimensionality and \eqref{eq:proof-gram-variation} yield
\[
\langle \mathcal G_\varepsilon z,z\rangle
\ge\frac c2\|z\|^2,
\qquad z\in \mathcal N_1,
\]
for some $c>0$ and all sufficiently small $|\varepsilon|$.

Relative to $R\oplus \mathcal N_1$ the remaining effective block is
\begin{equation}\label{eq:two-stage-schur}
\begin{pmatrix}
\varepsilon M_R-\varepsilon^2\mathcal G_{RR}(\varepsilon)
&-\varepsilon^2\mathcal G_{R\mathcal N_1}(\varepsilon)\\
-\varepsilon^2\mathcal G_{\mathcal N_1R}(\varepsilon)
&-\varepsilon^2\mathcal G_{\mathcal N_1\mathcal N_1}(\varepsilon)
\end{pmatrix}.
\end{equation}
Set $P_\varepsilon=\varepsilon M_R-
\varepsilon^2\mathcal G_{RR}(\varepsilon)$.  Since $M_R$ is
invertible and $\mathcal G_\varepsilon$ is uniformly bounded,
\[
P_\varepsilon=\varepsilon(M_R+O(\varepsilon)),
\qquad \|P_\varepsilon^{-1}\|=O(|\varepsilon|^{-1}).
\]
The homotopy
$\varepsilon M_R-t\varepsilon^2\mathcal G_{RR}(\varepsilon)$,
$0\le t\le1$, remains invertible for small nonzero
$\varepsilon$. Consequently,
$n_-(P_\varepsilon)=n_-(\varepsilon M_R)$.
Eliminating this block by a second Schur complement gives
\[
\begin{aligned}
T_\varepsilon={}&-\varepsilon^2\mathcal G_{\mathcal N_1\mathcal N_1}(\varepsilon)\\
&-\varepsilon^4\mathcal G_{\mathcal N_1R}(\varepsilon)
P_\varepsilon^{-1}\mathcal G_{R\mathcal N_1}(\varepsilon).
\end{aligned}
\]
The inverse bound and \eqref{eq:proof-gram-variation} imply
\[
T_\varepsilon
=-\varepsilon^2\bigl(\mathcal G_{0,\mathcal N_1}+O(|\varepsilon|)\bigr),
\]
where $\mathcal G_{0,\mathcal N_1}$ is positive definite.  In particular,
for $0<|\varepsilon|$ sufficiently small,
$T_\varepsilon$ is strictly negative on $\mathcal N_1$.
If $R=\{0\}$ or $\mathcal N_1=\{0\}$, the same argument applies
with the absent block omitted.  Both Schur reductions are
invertible congruences, so their inertias add:
\[
\begin{aligned}
\operatorname{ind}_-(\mathfrak d_{1+\varepsilon\eta})
&=n_-(\varepsilon M_R)+\dim \mathcal N_1\\
&=n_-(\varepsilon M_\eta)
+\operatorname{rank}(B_\eta|_{\ker M_\eta}).
\end{aligned}
\]
The complementary block is positive, and the only remaining
kernel is $\mathcal N_0$. Hence
$\operatorname{nul}(\mathfrak d_{1+\varepsilon\eta})
=\dim(\ker M_\eta\cap\ker B_\eta)$.
Finally, the decomposition of $E_1$ gives
$n_-(\varepsilon M_R)+\dim \mathcal N_1\le\dim E_1$.
This proves both identities in \eqref{thm:GH-index} and the
asserted upper bound.
\end{proof}

\begin{proposition}
\label{prop:quadratic-splitting}
Under the hypotheses of Theorem~\ref{thm:GH-index},
assume $M_{\eta}=0$. Let
$\nu_{1},\ldots,\nu_{d}\geq0$ be the eigenvalues,
counted with multiplicity, of
$\mathcal G_{0}=B_{\eta}^{*}C_{0}^{-1}B_{\eta}$
acting on the $L^{2}$-orthonormal first
eigenspace $E_{1}$. The generalized variational problem
\[
\mathfrak h_{1+\varepsilon\eta}(u,v)
=\kappa\,\mathfrak g_{1+\varepsilon\eta}(u,v),
\qquad v\in H^{2}(F)/\mathbb R,
\]
has exactly $d$ eigenvalues, counted with
multiplicity, converging to $\lambda_{1}$
as $\varepsilon\to0$. They can be numbered
so that
\begin{equation}\label{eq:quadratic-splitting}
\kappa_{j}(\varepsilon)
=\lambda_{1}
-\frac{\varepsilon^{2}}{\lambda_{1}}\nu_{j}
+O(|\varepsilon|^{3}),
\qquad j=1,\ldots,d.
\end{equation}
For every zero eigenvalue $\nu_{j}$, the
corresponding generalized eigenvalue is
exactly $\lambda_{1}$ for all sufficiently
small $\varepsilon$. Moreover, for every
sufficiently small nonzero $\varepsilon$,
exactly $\operatorname{rank}B_{\eta}$
eigenvalues lie strictly below $\lambda_{1}$.
\end{proposition}

\begin{proof}
At $\varepsilon=0$, the integrated Bochner identity gives
$\mathfrak h_1(u,v)=\langle Lu,Lv\rangle$ and
$\mathfrak g_1(u,v)=\langle Lu,v\rangle$.
Thus the generalized eigenvalues are the positive Laplace
eigenvalues on any closed connected manifold.
Thus $\lambda_1$ is isolated with multiplicity $d$.
Write $\kappa=\lambda_1+\delta$ and consider
\[
\mathcal Q_{\varepsilon,\delta}
=\mathfrak h_{1+\varepsilon\eta}
-(\lambda_1+\delta)\mathfrak g_{1+\varepsilon\eta}
=\mathfrak d_{1+\varepsilon\eta}
-\delta\mathfrak g_{1+\varepsilon\eta}.
\]
Relative to $E_1\oplus\mathcal V_+$, the assumption
$M_\eta=0$ and the orthogonality at $m$ give the blocks
\[
\begin{aligned}
\mathcal Q_{\varepsilon,\delta}|_{E_1\times E_1}
&=-\lambda_1\delta I+O(|\varepsilon\delta|),\\
\mathcal Q_{\varepsilon,\delta}|_{\mathcal V_+\times E_1}
&=\varepsilon B_\eta+O(|\varepsilon\delta|),\\
\mathcal Q_{\varepsilon,\delta}|_{\mathcal V_+\times\mathcal V_+}
&=C_0+O(|\varepsilon|+|\delta|).
\end{aligned}
\]
Coercivity of $C_0$ implies that the last block,
$D_{\varepsilon,\delta}$, is invertible for small parameters,
with
$D_{\varepsilon,\delta}^{-1}
=C_0^{-1}+O(|\varepsilon|+|\delta|)$.
A Schur reduction therefore gives the finite-dimensional
eigenvalue equation
\[
\det\mathcal S_{\varepsilon,\delta}=0,
\qquad
\mathcal S_{\varepsilon,\delta}
=-\lambda_1\delta I-\varepsilon^2\mathcal G_0
+O(|\varepsilon\delta|+|\varepsilon|^3),
\]
where $\mathcal G_0=B_\eta^*C_0^{-1}B_\eta$.
For a unit vector in the kernel of
$\mathcal S_{\varepsilon,\delta}$, this identity yields
$\lambda_1|\delta|\le C\varepsilon^2
+C|\varepsilon\delta|$, absorbing the last term gives
$\delta=O(\varepsilon^2)$. Hence
\[
\mathcal S_{\varepsilon,\delta}
=-\lambda_1\delta I-\varepsilon^2\mathcal G_0+O(|\varepsilon|^3).
\]
The fixed-space analytic realization in
Proposition~\ref{prop:spectral-realization}, followed by the
finite-dimensional min--max estimate, now gives
$-\lambda_1\delta_j(\varepsilon)
=\varepsilon^2\nu_j+O(|\varepsilon|^3)$,
which is \eqref{eq:quadratic-splitting}.

Finally, $\ker \mathcal G_0=\ker B_\eta$ by
\eqref{eq:gram-PW}. For $u$ in this kernel,
$\mathfrak d_{1+\varepsilon\eta}(u,v)=0$ for every $v$,
since both $M_\eta$ and $B_\eta u$ vanish.
Thus such vectors give exact generalized eigenvectors at
$\kappa=\lambda_1$. The remaining positive eigenvalues of
$\mathcal G_0$ yield strictly negative quadratic shifts, and
Theorem~\ref{thm:GH-index} gives their number as
$\operatorname{rank}B_\eta$.
\end{proof}

\begin{corollary}\label{cor:detect-criterion}
Assume that $F=G/H$ is homogeneous. For a fixed smooth real
mean-zero $\eta$, the inequality
$\mathfrak h_{1+\varepsilon\eta}\geq
\lambda_1\mathfrak g_{1+\varepsilon\eta}$ holds on all smooth
potentials for every sufficiently small nonzero $\varepsilon$
if and only if
\begin{equation}\label{eq:invisible-criterion}
M_\eta=0\quad\text{and}\quad B_\eta=0.
\end{equation}
If either of the two operators is nonzero, the inequality fails for
both signs of all sufficiently small nonzero admissible parameters.
If $\eta\in \mathscr E_\gamma$, condition \eqref{eq:invisible-criterion}
is equivalent to \eqref{eq:resonance-invisible}.
\end{corollary}
\begin{proof}
If $M_\eta\ne0$, Lemma~\ref{lem:tracefree} gives eigenvalues of
both signs, so $n_-(\varepsilon M_\eta)\geq1$ for either sign of
$\varepsilon$. If $M_\eta=0$ and $B_\eta\ne0$, the second term
in \eqref{thm:GH-index} is $\operatorname{rank}B_\eta>0$.
In both cases Theorem~\ref{thm:GH-index} gives a negative direction.
Conversely, if $M_\eta=B_\eta=0$, the Schur complement is zero
and the complementary block is positive for small $|\varepsilon|$,
so the full form is nonnegative. For a single Laplace mode,
Lemma~\ref{lem:first-row} gives the resonance characterization.
\end{proof}

\subsection{Spectral selection rules and rank}

We now return to the homogeneous setting $F=G/H$.
A permitted tensor-product type does not by itself yield a nonzero
matrix coefficient, and the spectral prefactor may kill a nonzero
product projection.  The Gram formula retains both effects.
When a density perturbation has symmetries, the finite-dimensional
operators admit a further representation-theoretic reduction.
Let
\[
G_\eta=\{g\in G:g\cdot\eta=\eta\}
\]
be its (closed) stabilizer in the isometric left $G$-action.  The
restriction $E_1^{\mathbb C}|_{G_\eta}$ decomposes as
\begin{equation}\label{eq:K-isotypic}
E_1^{\mathbb C}
=\bigoplus_{\sigma\in\widehat{G_\eta}}
U_\sigma\otimes\mathbb C^{m_\sigma}.
\end{equation}
The original Schur decomposition
$L|_{\mathcal W_\pi}=I_{V_\pi^*}\otimes A_\pi$ shows that
\begin{equation}\label{eq:K-mult}
m_\sigma=\sum_{\pi\in\widehat G}
\dim\ker(A_\pi-\lambda_1I)\,
\dim\operatorname{Hom}_{G_\eta}(U_\sigma,V_\pi^*|_{G_\eta}).
\end{equation}
Only finitely many $\pi$ contribute. These numbers are the actual
multiplicities in the first eigenspace, not fusion coefficients.

\begin{proposition}
\label{prop:multiplicity-index}
For a fixed smooth real mean-zero density perturbation $\eta$,
write the $G_\eta$-equivariant operators in \eqref{eq:MB-operators},
after complexification and the isotypic decomposition, as
\[
M_\eta=\bigoplus_\sigma I_{U_\sigma}\otimes a_\sigma,
\qquad
B_\eta=\bigoplus_\sigma I_{U_\sigma}\otimes b_\sigma,
\]
where $a_\sigma$ is a Hermitian matrix on the multiplicity space
and $b_\sigma$ maps that space to the corresponding multiplicity
space in $\mathcal V_+^*$.  Then, for sufficiently small
$\varepsilon\neq0$,
\begin{equation}\label{eq:K-index}
\displaystyle
\operatorname{ind}_-(\mathfrak d_{1+\varepsilon\eta})
=\sum_{\sigma\in\widehat{G_\eta}}(\dim U_\sigma)
\left[n_-(\varepsilon a_\sigma)
+\operatorname{rank}(b_\sigma|_{\ker a_\sigma})\right].
\end{equation}
The complexified index equals the original real index.  The
multiplicity matrices are computed from the actual triple
coefficients, not only from the list of allowed representation types.
\end{proposition}
\begin{proof}
The $G_\eta$-action preserves the metric, the density $\eta$, the
Laplacian, and thus the form $\mathfrak d_\eta$.  Consequently
$M_\eta$ and $B_\eta$ commute with $G_\eta$.  Schur's lemma yields the
block decompositions, including the corresponding isotypic
splitting of the target.  In each block, the kernel is
$U_\sigma\otimes\ker a_\sigma$ and the rank of the restricted
mixed map is $\dim U_\sigma$ times the rank of
$b_\sigma|_{\ker a_\sigma}$.  Theorem~\ref{thm:GH-index}
then gives \eqref{eq:K-index}.  Finally, a real symmetric form
and its complexification have the same negative index counted over
$\mathbb R$ and $\mathbb C$, respectively.
\end{proof}

The preceding formula depends on the ranks of $M_{\eta}$ and $B_{\eta}$. When $\eta$ has finite Peter-Weyl support, these ranks can be computed using finitely many $H$-fixed Clebsch-Gordan coefficients. This also applies when $\dim V_{\pi}^{H}>1$ and the operators $A_{\pi}$ have nontrivial spectral blocks.

\begin{proposition}\label{prop:finite-channel-test}
Let $\mathscr E\subset C^\infty(F,\mathbb R)\cap L_0^2(F,m)$
be a finite-dimensional space spanned by Laplace
eigenfunctions with finite Peter--Weyl expansions. Choose an eigenbasis $\{f_c\}_{c=1}^{s}$ of
$\mathscr E$, with $Lf_c=\gamma_cf_c$, and set
$\eta_t=\sum_{c=1}^{s}t_cf_c$. There exists a finite sum $Y$ of higher Laplace
eigenspaces containing $B_{\eta_t}(E_1)$ for every
$t\in\mathbb R^s$. Relative to orthonormal
eigenbases $\{u_i\}_{i=1}^{d}$ of $E_1$ and
$\{v_j\}_{j=1}^{r_Y}$ of $Y$, where $r_Y=\dim Y$, with
$Lv_j=\beta_jv_j$, define
\begin{align}
M(t)_{ij}
&=\frac14\sum_{c=1}^{s}t_c\gamma_c
(\gamma_c-2\lambda_1)
\int_F f_cu_ju_i\,dm,
\label{eq:finite-test-M}\\
B(t)_{ji}
&=\frac14\sum_{c=1}^{s}t_c
\bigl[(\beta_j-\gamma_c)^2-\lambda_1^2\bigr]
\int_F f_cu_iv_j\,dm.
\label{eq:finite-test-B}
\end{align}
Then, with
$R(t)=\bigl(\begin{smallmatrix}M(t)\\B(t)\end{smallmatrix}\bigr)$,
for every fixed $t$ and sufficiently small
nonzero admissible $\varepsilon$,
\begin{equation}\label{eq:finite-test-index}
\begin{aligned}
\operatorname{ind}_-(\mathfrak d_{1+\varepsilon\eta_t})
&=n_-(\varepsilon M(t))
+\operatorname{rank}R(t)-\operatorname{rank}M(t),\qquad
\operatorname{nul}(\mathfrak d_{1+\varepsilon\eta_t}) =d-\operatorname{rank}R(t).
\end{aligned}
\end{equation}
The entries are computable from finitely many
$H$-fixed Clebsch--Gordan coefficients.
Moreover, the parameter space admits a
semialgebraic stratification on which the
ranks and inertia appearing above are constant.
\end{proposition}

\begin{proof}
For fixed $f_c$ and $u_i$, the product $f_cu_i$ is a finite sum
of matrix coefficients belonging to the finite tensor products of
their $G$-types.  Thus the product has finite Peter--Weyl support.
Each isotypical component splits under $L$ into finitely many
spectral subblocks, because its $H$-fixed multiplicity space is
finite-dimensional.  We may therefore take $Y$ to be the sum of all
spectral subblocks of eigenvalue greater than $\lambda_1$ that
occur in the finitely many products $f_cu_i$.  Equation
\eqref{eq:matrix-Bbeta} now shows that the functional
$B_{\eta_t}u_i$ is represented by an $L^2$ function in $Y$,
for every $t$, components killed by the resonance factor may
optionally be omitted from $Y$.  Equations
\eqref{eq:matrix-M} and \eqref{eq:matrix-Bbeta} give
\eqref{eq:finite-test-M} and \eqref{eq:finite-test-B}.
The coefficients in these equations are triple products of
matrix coefficients, hence are computed by Lemma~\ref{lem:CG-explicit},
including all multiplicity directions.

For two finite-dimensional linear maps
$M:E_1\to E_1$ and $B:E_1\to Y$, consider the combined map
$R=(M,B):E_1\to E_1\oplus Y$. Projection from
$\operatorname{im}R$ onto $\operatorname{im}M$ has kernel
isomorphic to $B(\ker M)$. Therefore
\[
\operatorname{rank}R
=\operatorname{rank}M+\operatorname{rank}(B|_{\ker M}).
\]
Likewise $\ker R=\ker M\cap\ker B$.  Substitution into
Theorem~\ref{thm:GH-index} proves \eqref{eq:finite-test-index}.
This identity is valid even if $M$ is singular and its kernel
changes with $t$, no continuously chosen kernel basis is required.

Every matrix entry is linear in the real coordinates $t_c$.
A rank condition is thus a condition on polynomial minors, and
signature is constant on each connected part of a fixed-rank
stratum because the nonzero eigenvalues of a real symmetric
matrix cannot change sign without crossing zero.  The ranks
are maximal wherever a selected nonzero maximal minor does not
vanish, the exceptional set is a proper real algebraic subset.
A finite Boolean combination of signs and zero sets of minors and
characteristic-polynomial coefficients gives the asserted
semialgebraic decomposition.  Positivity of $1+\varepsilon\eta_t$
is ensured for each fixed $t$ by taking $|\varepsilon|$ small.
\end{proof}

\begin{example}\label{ex:resonant-perturbation}
The vanishing of the spectral obstructions does not
necessarily force the perturbation to be constant.
Indeed, consider $F=S^{1}\times S^{2}$ with the
unit product metric, for which $\lambda_{1}=1$
and $E_{1}$ consists of the first eigenfunctions
on $S^{1}$. Let $\eta$ be a nonzero degree-one
spherical harmonic on $S^{2}$, extended
constantly along $S^{1}$. Then $L\eta=2\eta$,
and for every $u\in E_{1}$, the product
$\eta u$ belongs to $\mathscr E_{3}$.

This interaction is resonant: with
$\gamma=2$ and $\beta=3$, the coefficients
in Lemma~\ref{lem:first-row} satisfy
\[
\begin{aligned}
\gamma(\gamma-2\lambda_{1})
&=2(2-2)=0,\\
(\beta-\gamma)^{2}-\lambda_{1}^{2}
&=(3-2)^{2}-1=0.
\end{aligned}
\]
Consequently, $M_{\eta}=B_{\eta}=0$.
The Schur complement therefore vanishes
identically, while the higher-mode block
$C_{\varepsilon}$ remains positive for
sufficiently small $|\varepsilon|$. Hence
\[
\mathfrak h_{1+\varepsilon\eta}(u,u)
\geq\lambda_{1}
\mathfrak g_{1+\varepsilon\eta}(u,u)
\]
for every admissible $u$, with equality
along $E_{1}$, even though
$(1+\varepsilon\eta)m$ is non-homogeneous
whenever $\varepsilon\neq0$.

Thus the preservation of the Hessian bound
is not restricted to homogeneous measures.
It can also arise from spectral resonance,
when the perturbation couples the first
eigenspace only to modes for which the
obstruction coefficients vanish. This
illustrates why the loss of
the fixed $\lambda_{1}$-Hessian inequality is governed by
the actual spectral couplings, rather than
by the nonconstancy of the density alone.
\end{example}

\section{Rigidity and the local index on round spheres}
\label{sec:spheres}

Consider the unit round sphere
\[
S^d=SO(d+1)/SO(d),\qquad d\geq2,
\]
with normalized volume. Let $\mathscr H_n$ be the real spherical
harmonics of degree $n$, with complexification understood when
using representation coefficients. Then
\[
L|_{\mathscr H_n}=\lambda_n I,\qquad
\lambda_n=n(n+d-1),\qquad E_1=\mathscr H_1,\qquad \lambda_1=d.
\]
For $\eta\in C^\infty(S^d,\mathbb R)$ of mean zero, let $\eta_2$
be its orthogonal projection onto $\mathscr H_2$, and write
\[
\eta_2(x)=x^{\mathsf T}Q_\eta x,\qquad
Q_\eta=Q_\eta^{\mathsf T},\quad \operatorname{tr}Q_\eta=0.
\]
This identifies the degree-two component uniquely, $Q_\eta$ may
vanish for a nonzero perturbation.

The geometry of first harmonics gives a useful alternative to
computing every coupling coefficient. Since $\Hess u=-ug$ for
$u\in E_1$ and $\Ric=(d-1)g$, integration by parts yields
\begin{equation}\label{eq:sphere-density-gradient}
\mathfrak d_f(u,v)
=\int_{S^d}u\,\langle\nabla f,\nabla v\rangle\,dm,
\qquad u\in E_1,\quad v\in C^\infty(S^d,\mathbb R).
\end{equation}
Here $f$ is any smooth real function and $\mathfrak d_f$ is the
linear extension of the defect in the density variable. Indeed,
its integrand before integration by parts is
$f(-u\Delta v-\langle\nabla u,\nabla v\rangle)$.
For complex potentials the identity is extended sesquilinearly.

\begin{lemma}\label{lem:sphere-one-direction}
Let $f\in C^\infty(S^d,\mathbb R)$. If there is a nonzero
$u\in E_1$ such that $\mathfrak d_f(u,v)=0$ for every smooth
real $v$, then $f$ is constant. Consequently, for every nonzero
smooth real mean-zero $\eta$,
\[
\ker M_\eta\cap\ker B_\eta=\{0\}.
\]
\end{lemma}
\begin{proof}
Equation~\eqref{eq:sphere-density-gradient} implies
$\operatorname{div}(u\nabla f)=0$ in distributions, hence
pointwise. After a rotation and a nonzero scaling take
$u=x_{d+1}$. On the closed hemisphere
$\overline{\Omega_+}=\{x_{d+1}\geq0\}$, the divergence theorem gives
\[
0=\int_{\Omega_+}f\,\operatorname{div}(u\nabla f)\,dm
=-\int_{\Omega_+}u|\nabla f|^2\,dm,
\]
because the boundary term $fu\,\partial_\nu f$ vanishes on the
equator. Since $u>0$ in $\Omega_+$, $f$ is constant there.
Applying the same argument on $\Omega_-=\{u<0\}$ and using
continuity across the equator shows that the constants agree.
For the last assertion, a nonzero vector in
$\ker M_\eta\cap\ker B_\eta$ would pair trivially with
$E_1\oplus\mathcal V_+$ and with constants, forcing $\eta$ to be
constant and hence zero.
\end{proof}

\begin{theorem}\label{thm:sphere-index}
On the unit round $S^d$, $d\ge2$, the inequality
$\mathfrak h_\rho\ge d\,\mathfrak g_\rho$ holds for a smooth
positive probability density $\rho$ if and only if $\rho\equiv1$.
Moreover, for every nonzero smooth real mean-zero $\eta$ and every
sufficiently small nonzero $\varepsilon$ with
$1+\varepsilon\eta>0$,
\begin{equation}\label{eq:sphere-index}
\begin{aligned}
\operatorname{ind}_-(\mathfrak d_{1+\varepsilon\eta})
&=n_-(\varepsilon Q_\eta)+\dim\ker Q_\eta,\\
\operatorname{nul}(\mathfrak d_{1+\varepsilon\eta})&=0.
\end{aligned}
\end{equation}
\end{theorem}
\begin{proof}
At $\rho=1$ the asserted inequality follows from the homogeneous
spectral decomposition.
Conversely, if $\mathfrak d_\rho$ is nonnegative, its restriction
to $E_1$ is nonnegative and trace-free by
Lemma~\ref{lem:tracefree}, so this restriction vanishes.
For $u\in E_1$ and any smooth $v$, nonnegativity of
$\mathfrak d_\rho(u+tv,u+tv)$ for all real $t$ then forces
$\mathfrak d_\rho(u,v)=0$. Lemma~\ref{lem:sphere-one-direction}
makes $\rho$ constant, and normalization gives $\rho=1$.

For the local statement, choose the orthonormal real basis
$u_i=\sqrt{d+1}\,x_i$ of $E_1$, $1\leq i\leq d+1$.
The pointwise defect integrand is
\[
(d+1)^2x_ix_j-(d+1)\delta_{ij}.
\]
Only degrees zero and two occur, so mean zero implies
$M_\eta=M_{\eta_2}$. The spherical fourth moments
\[
\int_{S^d}x_ix_jx_kx_\ell\,dm
=\frac{\delta_{ij}\delta_{k\ell}
+\delta_{ik}\delta_{j\ell}+\delta_{i\ell}\delta_{jk}}
{(d+1)(d+3)}
\]
give
\begin{equation}\label{eq:sphere-moment-matrix}
M_\eta=\frac{2(d+1)}{d+3}Q_\eta.
\end{equation}
Lemma~\ref{lem:sphere-one-direction} says that $B_\eta$ is
injective on $\ker M_\eta$. Theorem~\ref{thm:GH-index}
therefore gives both identities in \eqref{eq:sphere-index}.
No finite-support assumption on $\eta$ has been used.
\end{proof}

The coupling characterization implicit in the rigidity argument
is useful independently of the local index formula.

\begin{corollary}\label{cor:sphere-critical-row}
On the unit round $S^d$, $d\ge2$, a smooth positive probability
density $\rho$ satisfies $\rho\equiv1$ if and only if
$\mathfrak d_\rho(E_1,\mathscr H_k)=0$ for every $k\ge1$.
\end{corollary}
\begin{proof}
At $\rho=1$ the critical row vanishes by the homogeneous spectral
decomposition. Conversely, the stated vanishing, together with
convergence of spherical harmonic expansions in $H^2$ and continuity
of the defect form, gives $\mathfrak d_\rho(u,v)=0$ for every
$u\in E_1$ and every smooth real $v$, constants pair trivially.
Lemma~\ref{lem:sphere-one-direction} forces $\rho$ to be constant,
and normalization gives $\rho=1$.
\end{proof}

For a fixed nonzero perturbation $\eta$, the local index in
Theorem~\ref{thm:sphere-index} depends only on its degree-two
component. The size of the neighborhood in $\varepsilon$ where
the formula holds may nevertheless depend on the full perturbation.
The next corollary records the resulting counts for individual
harmonic degrees and for perturbations with no degree-two component.

\begin{corollary}\label{cor:sphere-harmonic-index}
If $0\ne\eta\in\mathscr H_n$ is real and $n\geq1$, then for
sufficiently small nonzero $\varepsilon$,
\[
\operatorname{ind}_-(\mathfrak d_{1+\varepsilon\eta})
=\begin{cases}
d+1,&n\ne2,\\
n_-(\varepsilon Q_\eta)+\dim\ker Q_\eta,&n=2.
\end{cases}
\]
More generally, every nonzero smooth mean-zero $\eta$ with
$\eta_2=0$ has local index $d+1$. If $\eta_2\ne0$, the index
belongs to $\{1,\ldots,d\}$.
\end{corollary}
\begin{proof}
Apply \eqref{eq:sphere-index}. If $Q_\eta=0$, its kernel has
dimension $d+1$. Otherwise its trace-free spectrum contains both
positive and negative eigenvalues.
\end{proof}

The higher modes can therefore change the sizes of the eigenvalue
shifts without changing this local count for a fixed nonzero
quadratic component. When that component vanishes, all $d+1$
critical eigenvalues move below $d$ at second order, with positive
Gram eigenvalues by Lemma~\ref{lem:sphere-one-direction} and
Proposition~\ref{prop:quadratic-splitting}.
This local assertion is separate from the global rigidity above.

\begin{example}\label{ex:sphere-mixed-degrees}
On $S^2$, take $\eta(x,y,z)=x^2-y^2+z$. Its degree-two matrix is
$Q_\eta=\operatorname{diag}(1,-1,0)$, while $z$ is a degree-one
component. Formula~\eqref{eq:sphere-index} gives index two and
nullity zero for both signs of sufficiently small nonzero
$\varepsilon$: one negative direction comes from the first-order
matrix, and the other from its one-dimensional kernel through
the mixed coupling. The formula applies equally when further
smooth components of degrees other than zero and two are added.
\end{example}

\subsection{Harmonic channels and a quantitative example}

The index proof does not require individual channel calculations.
For comparison with the general representation formulas, the
spherical products identify those channels explicitly
\cite[Ch.~5]{AxlerBourdonRamey2001}.

\begin{lemma}\label{lem:sphere-channels}
Write $\rho=1+\sum_{k\geq1}\eta_k$ with
$\eta_k\in\mathscr H_k$. For $n\geq1$, $u\in\mathscr H_1$, and
$v\in\mathscr H_{n+1}$, the sesquilinear defect satisfies
\begin{equation}\label{eq:sphere-channels}
\begin{aligned}
\mathfrak d_\rho(u,v)
={}&n(n+d)\int_{S^d}\eta_nu\overline v\,dm\\
&+(n+1)(n+d+1)\int_{S^d}\eta_{n+2}u\overline v\,dm.
\end{aligned}
\end{equation}
The two maps from $\mathscr H_n$ and $\mathscr H_{n+2}$ to these
sesquilinear pairings are injective and have disjoint images.
\end{lemma}
\begin{proof}
For a harmonic homogeneous polynomial $P_m$ on $\mathbb R^{d+1}$,
the Fischer decomposition is
\begin{equation}\label{eq:spherical-product-explicit}
x_iP_m=U_{m+1}+\frac{|x|^2}{2m+d-1}\partial_iP_m,
\end{equation}
where $U_{m+1}$ is harmonic. Thus $u\overline v$ contains only
degrees $n$ and $n+2$. Substitution of
$\lambda_k=k(k+d-1)$ into Lemma~\ref{lem:first-row} gives the
coefficients in \eqref{eq:sphere-channels}.

Both pairing maps are $SO(d+1)$-equivariant. The lower-degree
projection in \eqref{eq:spherical-product-explicit} is nonzero
for some $i$, since a nonconstant polynomial has a nonzero
partial derivative. The upper-degree projection is also nonzero:
if $U_{m+1}=0$ for a nonzero $P_m$ and some $x_i$, then
$|x|^2$ divides $x_iP_m$. In $d+1\geq3$ variables this quadratic
is irreducible over $\mathbb C$, so it divides $P_m$.
A nonzero harmonic polynomial cannot have that factor, since
multiplication by $|x|^2$ is adjoint to the Laplacian for the
Fischer inner product. Indeed $P_m=|x|^2R$ would imply
$\|P_m\|_{\mathrm F}^2=\langle\Delta P_m,R\rangle_{\mathrm F}=0$.
Both maps are therefore nonzero. The irreducibility and
inequivalence of $\mathscr H_n$ and $\mathscr H_{n+2}$ give
injectivity and disjoint images.
\end{proof}

\begin{example}
\label{ex:sphere2-gram}
Consider $\eta(x,y,z)=z$ on $S^2$, where $\gamma=\lambda_1=2$. Since the product of two first harmonics contains only degrees $0$ and $2$, while $\eta$ has degree $1$, the first-order matrix $M_\eta$ vanishes. To compute $B_\eta$, consider the normalized basis of the first eigenspace given by
$$
\{\sqrt{3}\,x,\sqrt{3}\,y,\sqrt{3}\,z\}
$$
Multiplication by $z$ produces both a constant component and a degree-two component, but only the latter contributes to $B_\eta$, with coefficient $\tfrac14[(6-2)^2-2^2]=3$. To compute the corresponding second-order contribution, we also use the homogeneous defect on $\mathscr H_2$, which equals $6(6-2)=24$. 
Thus,
\[
B_{\eta,6}u_x=3\sqrt{3}\,xz,\qquad
B_{\eta,6}u_y=3\sqrt{3}\,yz,\qquad
B_{\eta,6}u_z=3\sqrt{3}\left(z^2-\frac{1}{3}\right).
\]
Using the spherical moments
$\int x^2z^2\,dm=1/15$ and
$\int z^4\,dm=1/5$, we obtain the diagonal entries explicitly:
\[
\frac{\|3\sqrt{3}\,xz\|_{L^2}^2}{24}
=\frac{27}{24}\cdot\frac{1}{15}=\frac{3}{40},
\qquad
\frac{\|3\sqrt{3}(z^2-\frac{1}{3})\|_{L^2}^2}{24}
=\frac{27}{24}\cdot\frac{4}{45}=\frac{1}{10}.
\]
We therefore obtain the exact spectral Gram operator:
\[
B_\eta^*C_0^{-1}B_\eta
=\operatorname{diag}\left(
\frac{3}{40},\frac{3}{40},\frac{1}{10}
\right),
\]
which is strictly positive. By the general theorem, this gives
$\operatorname{ind}_-(\mathfrak d_{1+\varepsilon z})=3$
for all sufficiently small nonzero $\varepsilon$, with all three
negative directions arising at second order, even though the
first-order block vanishes identically. Moreover,
Proposition~\ref{prop:quadratic-splitting} gives the generalized
entropy-Hessian eigenvalues near $2$:
\[
\kappa_x(\varepsilon)=\kappa_y(\varepsilon)
=2-\frac{3}{80}\varepsilon^2+O(|\varepsilon|^3),
\qquad
\kappa_z(\varepsilon)
=2-\frac{1}{20}\varepsilon^2+O(|\varepsilon|^3).
\]
These coefficients are the Gram eigenvalues divided by $\lambda_1=2$, not the eigenvalues of the defect form itself.
\end{example}

\section{Complete flags: an index calculation with multiplicities}
\label{sec:complete-flags}

Complete flags exhibit a feature that is invisible in the basic
Fourier examples: one irreducible representation may contribute
several independent copies to the same Laplace eigenspace. For
$N\ge3$, the complete complex flag manifold is
\begin{equation}\label{eq:flag-definition-reader}
\Fl_N=\bigl\{(V_1,\ldots,V_{N-1}):
0\subset V_1\subset\cdots\subset V_{N-1}\subset\mathbb C^N,
\ \dim_{\mathbb C}V_j=j\bigr\}.
\end{equation}
The action of $SU(N)$ on nested subspaces is transitive, and the
stabilizer of the standard coordinate flag is its diagonal maximal
torus $T=T^{N-1}$. Indeed, the columns of a unitary matrix determine
an adapted orthonormal basis, while changing their phases preserves
every subspace of the flag. Thus $\Fl_N=SU(N)/T$ and
$\dim_{\mathbb R}\Fl_N=N(N-1)$.

We equip $SU(N)$ with the bi-invariant inner product
$Q(X,Y)=-\operatorname{tr}(XY)$ on $\mathfrak{su}(N)$ and
$\Fl_N$ with the induced normal homogeneous metric. The quotient
map is a Riemannian submersion with totally geodesic fibers, hence
the operator $L=-\Delta$ on the quotient is the quadratic Casimir restricted
to right-$T$-invariant functions. A representation contributes to
$L^2(\Fl_N)$ with multiplicity equal to the dimension of its
$T$-fixed subspace. For the adjoint representation this is the
Cartan space of traceless diagonal matrices, of dimension $N-1$.
The first positive eigenvalue is $2N$, and we verify below, using
the highest-root order and a concrete Casimir normalization, that
its full eigenspace consists precisely of these adjoint copies.
The resulting multiplicities are indispensable to the index count,
a calculation on one adjoint copy would miss most of the directions.

Choose the traceless diagonal matrix
$A=\operatorname{diag}(1,\ldots,1,-(N-1))$. Since it commutes
with $T$, the matrix coefficient $\eta_N(gT)=\operatorname{tr}
(AgAg^{-1})$ is well defined. Writing $A=I-NP_N$, where $P_N$
is the orthogonal projection onto the last coordinate axis, gives
\begin{equation}\label{eq:flag-eta}
\eta_N(gT)=N^2|g_{NN}|^2-N.
\end{equation}
Haar invariance yields $\int_{\Fl_N}|g_{NN}|^2\,dm_N=1/N$,
so $\int\eta_N\,dm_N=0$, moreover,
$-N\le\eta_N\le N(N-1)$. Therefore the densities
$\rho_{\varepsilon,N}=1+\varepsilon\eta_N$ are positive
for $-1/[N(N-1)]<\varepsilon<1/N$.
Unlike a general perturbation, this direction is accessible to
exact calculation: its cubic moment reduces to moments of one Haar
matrix entry, and the resulting defect factors across the adjoint
representation and its Cartan multiplicity space. We carry out
both calculations explicitly after identifying the first eigenspace
\cite{Yamaguchi1979,Haber2021}.

\begin{theorem}\label{thm:flag-index}
For every $N\ge3$, the normal complete flag $\Fl_N$ satisfies
\[
\lambda_1(\Fl_N)=2N,
\qquad \dim E_1=(N^2-1)(N-1).
\]
For each fixed $N$ and all sufficiently small admissible
$\varepsilon\ne0$, the negative index of the entropy-Hessian defect
is
\begin{equation}\label{eq:flag-index}
\operatorname{ind}_-(\mathfrak d_{\rho_{\varepsilon,N}})
=\begin{cases}
N(N-2)^2+2N-1,&\varepsilon>0,\\
(N-2)(3N-1),&\varepsilon<0.
\end{cases}
\end{equation}
\end{theorem}

Here admissibility means $1+\varepsilon\eta_N>0$, and the
smallness threshold may depend on $N$. We give a proof in three
steps: identification of the first eigenspace, an exact Haar
triple-product identity, and diagonalization of its defect matrix
with all representation multiplicities retained.

Let $\mathfrak a_N$ denote the real vector space of
traceless Hermitian $N\times N$ matrices, equipped with
the inner product $\langle A,B\rangle=\operatorname{tr}(AB)$,
and let $\mathfrak t_N$ be its diagonal subspace.
For $A\in\mathfrak a_N$ and $D\in\mathfrak t_N$,
consider the matrix coefficient
\[
f_{A,D}(gT)=\operatorname{tr}(AgDg^{-1}).
\]
Since $D$ commutes with every element of $T$, this
expression is independent of the representative $g$
and therefore defines a smooth real-valued function
on $\Fl_N$.

\begin{lemma}
\label{lem:flag-first-eigen}
The first positive eigenvalue of $L=-\Delta$ on
$\Fl_N$ is $2N$, and its eigenspace $E_1$ is naturally
isomorphic to $\mathfrak a_N\otimes\mathfrak t_N$
through the correspondence
\[
A\otimes D\longmapsto f_{A,D}.
\]
Moreover, for any $A,B\in\mathfrak a_N$ and
$D,E\in\mathfrak t_N$, the orthogonality relation
\begin{equation}\label{eq:flag-schur-orth}
\int_{\Fl_N}f_{A,D}f_{B,E}\,dm_N
=\frac{\operatorname{tr}(AB)\operatorname{tr}(DE)}
{N^2-1}
\end{equation}
holds. In particular,
$\dim E_1=(N^2-1)(N-1)$.
\end{lemma}

\begin{proof}
For the normal metric induced by
$Q(X,Y)=-\operatorname{tr}(XY)$, the nonnegative operator
$L=-\Delta$ on $\Fl_N=SU(N)/T$ is the positive quadratic Casimir.
Its spectral decomposition is therefore determined
by the irreducible representations of $SU(N)$
containing a nonzero $T$-fixed vector
\cite{Yamaguchi1979,FultonHarris1991,BrockerDieck1985}.

Let $\{X_a\}_{a=1}^{N^2-1}$ be a $Q$-orthonormal
basis of $\mathfrak{su}(N)$ and put $T_a=-iX_a$.
The matrices $T_a$, together with $I/\sqrt N$,
form an orthonormal basis of the Hermitian matrices
for the trace pairing. Complexifying the corresponding
completeness relation gives
\[
\sum_a(T_a)_{ij}(T_a)_{k\ell}
=\delta_{i\ell}\delta_{jk}
-\frac1N\delta_{ij}\delta_{k\ell}.
\]
Multiplying by $Z_{jk}$ and summing over $j,k$,
we obtain
\[
\sum_a T_aZT_a
=\operatorname{tr}(Z)I-\frac1N Z.
\]
Since $X_a=iT_a$, it follows that
\[
\sum_aX_aZX_a
=\frac1N Z-\operatorname{tr}(Z)I,
\qquad
\sum_aX_a^2=-\frac{N^2-1}{N}I.
\]
Consequently, for every traceless Hermitian matrix $Y$,
\begin{align*}
-\sum_a[X_a,[X_a,Y]]
&=-\sum_a(X_a^2Y-2X_aYX_a+YX_a^2)\\
&=\frac{2(N^2-1)}{N}Y+\frac2N Y\\
&=2NY.
\end{align*}

Thus, the Casimir operator acts on the adjoint representation
by $2N$. It remains to show that this is the smallest positive
Casimir eigenvalue among representations containing
weight zero. Let
$\lambda=(\lambda_1,\ldots,\lambda_N)$ be a nonzero
dominant highest weight in the root lattice, so that
$\lambda_1\ge\cdots\ge\lambda_N$ are integers and
$\sum_i\lambda_i=0$.
Its partial sums $S_j=\sum_{i=1}^j\lambda_i$
satisfy $S_j\ge1$ for $1\le j<N$: by dominance
they are nonnegative, and equality would force
all entries of $\lambda$ to vanish.

Writing $\vartheta=(1,0,\ldots,0,-1)$ for the
highest root, we have
\[
\lambda-\vartheta
=\sum_{j=1}^{N-1}(S_j-1)\alpha_j,
\qquad
\alpha_j=e_j-e_{j+1}.
\]
The Casimir eigenvalue associated with $\lambda$ is
$c(\lambda)=\langle\lambda,\lambda+2\varrho\rangle$,
where $\varrho$ is the half-sum of positive roots.
Since $\langle\alpha_j,\varrho\rangle=1$,
\begin{align*}
c(\lambda)-c(\vartheta)
&=\langle\lambda-\vartheta,
          \lambda+\vartheta+2\varrho\rangle\\
&=\sum_{j=1}^{N-1}(S_j-1)
  \bigl(\lambda_j-\lambda_{j+1}
  +\vartheta_j-\vartheta_{j+1}+2\bigr).
\end{align*}
Every term is nonnegative, and the sum vanishes
only when $\lambda=\vartheta$. Moreover,
\[
c(\vartheta)
=\langle\vartheta,\vartheta\rangle
 +2\langle\vartheta,\varrho\rangle
=2+2(N-1)=2N.
\]
Hence the adjoint representation gives precisely
the first positive eigenspace.

Its $T$-fixed subspace is $\mathfrak t_N$,
of dimension $N-1$. The corresponding matrix
coefficients are $f_{A,D}$, with
$A\in\mathfrak a_N$ and $D\in\mathfrak t_N$.

Finally, the adjoint representation on
$\mathfrak a_N$ is real irreducible.
For $D,E\in\mathfrak t_N$, Schur's lemma gives
\[
\int_{SU(N)}
(\operatorname{Ad}_g D)\otimes
(\operatorname{Ad}_g E)\,dg
=\frac{\operatorname{tr}(DE)}{N^2-1}
 \operatorname{Id}_{\mathfrak a_N},
\]
where $u\otimes v$ denotes the operator
$Z\mapsto\langle Z,v\rangle u$.
Pairing this identity with $A,B\in\mathfrak a_N$
yields
\[
\int_{\Fl_N}f_{A,D}f_{B,E}\,dm_N
=\frac{\operatorname{tr}(AB)\operatorname{tr}(DE)}
       {N^2-1}.
\]
In particular, the matrix-coefficient map
$\mathfrak a_N\otimes\mathfrak t_N\to E_1$
is injective and, by the Peter--Weyl decomposition,
surjective. Therefore
$\dim E_1=(N^2-1)(N-1)$.
\end{proof}

\begin{lemma}
\label{lem:flag-cubic}
For $B_1,B_2,B_3\in\mathfrak a_N$ and
$D_1,D_2,D_3\in\mathfrak t_N$, define
\[
d(X,Y,Z)=\operatorname{tr}\bigl(X(YZ+ZY)\bigr).
\]
Then
\begin{equation}\label{eq:flag-cubic}
\int_{\Fl_N}\prod_{j=1}^3f_{B_j,D_j}\,dm_N
=c_N d(B_1,B_2,B_3)d(D_1,D_2,D_3),
\qquad
c_N=\frac{N}{2(N^2-1)(N^2-4)}.
\end{equation}
\end{lemma}
\begin{proof}

Identify the adjoint representation with the real orthogonal
representation on $\mathfrak a_N$ given by conjugation, so that
Haar integration of its third tensor power gives the orthogonal
projection onto the invariant tensors. To describe these tensors,
we complexify $\mathfrak a_N$ to $\mathfrak{sl}_N(\mathbb C)$
and consider conjugation-invariant trilinear forms on traceless
matrices $X,Y,Z$. By the first fundamental theorem for simultaneous
conjugation in degree three \cite{Procesi1976}, equivalently
the permutation-contraction description from Schur--Weyl duality,
every such form is a linear combination of products of traces
of words. Since $\operatorname{tr}X=\operatorname{tr}Y=
\operatorname{tr}Z=0$, only the cyclic words
$\operatorname{tr}(XYZ)$ and $\operatorname{tr}(XZY)$ remain,
and they are independent for $N\ge3$. Their sum gives
$d(X,Y,Z)=\operatorname{tr}(X(YZ+ZY))$, while their
difference is proportional to
$f(X,Y,Z)=i\operatorname{tr}(X[Y,Z])$.
The first is symmetric and the second alternating, so they
are orthogonal under the natural inner product on tensor
powers, since interchanging two factors preserves $d$
and changes the sign of $f$.

To obtain the Haar integral from this classification, we regard
$d$ and $f$ as elements of the real Hilbert tensor space
$\mathfrak a_N^{\otimes 3}$, using an orthonormal basis of
$\mathfrak a_N$. The Haar average
$P=\int_{SU(N)}\operatorname{Ad}_g^{\otimes 3}\,dg$
is an orthogonal projection, since invariance of Haar measure
gives $P^2=P$, while invariance under inversion gives
$P^*=P$. As its range is precisely the space of invariant
tensors, we obtain
$$
P=\frac{d\otimes d}{\|d\|^2}
  +\frac{f\otimes f}{\|f\|^2},
$$
where the rank-one operators act by
$X\mapsto\langle X,d\rangle d/\|d\|^2$,
and similarly for $f$.

Applying this projector to
$D_1\otimes D_2\otimes D_3$ and pairing against
$B_1\otimes B_2\otimes B_3$ gives
\[
 \int_{SU(N)}\prod_{j=1}^3
  \langle\operatorname{Ad}_gD_j,B_j\rangle\,dg
 =\frac{d(B_1,B_2,B_3)d(D_1,D_2,D_3)}{\|d\|^2}
  +\frac{f(B_1,B_2,B_3)f(D_1,D_2,D_3)}{\|f\|^2}.
\]
Because $D_1,D_2,D_3$ are diagonal, they commute, and therefore
$f(D_1,D_2,D_3)=0$.  Thus the integral is precisely
$c_N d(B_1,B_2,B_3)d(D_1,D_2,D_3)$ with
$c_N=\|d\|^{-2}>0$.  The remaining task is to calculate
this numerical constant in the trace normalization used here.

To determine the constant, let $g$ be Haar-distributed on $SU(N)$
and put $X=|g_{NN}|^2$. Its last column is uniform on the unit
sphere in $\mathbb C^N$, equivalently a normalized complex
Gaussian vector. Thus $X$ has the Beta$(1,N-1)$ distribution and
\[
\mathbb E[X^k]=(N-1)\int_0^1x^k(1-x)^{N-2}\,dx
=\frac{k!(N-1)!}{(N+k-1)!}.
\]
These expectations equal the corresponding integrals on $\Fl_N$,
since $|g_{NN}|^2$ is right $T$-invariant. In particular,
\[
\mathbb E[X]=\frac1N,\qquad
\mathbb E[X^2]=\frac{2}{N(N+1)},\qquad
\mathbb E[X^3]=\frac{6}{N(N+1)(N+2)}.
\]
By \eqref{eq:flag-eta}, $f_{A,A}=N^2X-N$.  Because
integration over $\Fl_N$ agrees with the Haar expectation,
the integral of $f_{A,A}^3$ is $\mathbb E[(N^2X-N)^3]$.
Expanding the cube and using linearity of expectation gives
\begin{align*}
\int_{\Fl_N}f_{A,A}^3\,dm_N
&=N^6\mathbb E[X^3]-3N^5\mathbb E[X^2]
+3N^4\mathbb E[X]-N^3\\
&=\frac{2N^3(N-1)(N-2)}{(N+1)(N+2)}.
\end{align*}
On the other hand,
\[
d(A,A,A)=2\operatorname{tr}(A^3)
=-2N(N-1)(N-2).
\]
Dividing the cubic moment by $d(A,A,A)^2$ yields
\[
c_N=\frac{N}{2(N-1)(N-2)(N+1)(N+2)}
=\frac{N}{2(N^2-1)(N^2-4)},
\]
as claimed.
\end{proof}

The cubic moment is now normalized exactly, so the restriction of
the first variation to the critical eigenspace can be expressed
as a tensor product of concrete operators. This is the point at
which representation-theoretic multiplicity becomes a computable
signature, as the proof of Theorem~\ref{thm:flag-index} shows.

We now prove Theorem~\ref{thm:flag-index}.

\begin{proof}
Let $T_1,\ldots,T_{N^2-1}$ be an orthonormal real basis of
$\mathfrak a_N$ and $D_1,\ldots,D_{N-1}$ an orthonormal basis of
$\mathfrak t_N$.  Lemma~\ref{lem:flag-first-eigen} gives the
orthonormal basis
\[
\Psi_{a,i}=\sqrt{N^2-1}\,f_{T_a,D_i},
\qquad 1\le a\le N^2-1,\quad1\le i\le N-1,
\]
of the first eigenspace. Since $\eta_N\in E_1$ and
$\lambda_1=2N$, the first-eigenspace specialization of the
coupling formula \eqref{eq:D-triple} yields
\begin{equation}\label{eq:flag-first-defect}
\mathfrak d_{1+\varepsilon\eta_N}(u,v)
=-\varepsilon N^2\int_{\Fl_N}\eta_N uv\,dm_N,
\qquad u,v\in E_1.
\end{equation}
We calculate one matrix entry before using tensor notation.
The normalized basis is $\Psi_{a,i}=\sqrt{N^2-1}f_{T_a,D_i}$,
so Lemma~\ref{lem:flag-cubic} gives
\begin{align*}
\langle M_{\eta_N}\Psi_{a,i},\Psi_{b,j}\rangle
&=-N^2(N^2-1)
  \int_{\Fl_N}f_{A,A}f_{T_a,D_i}f_{T_b,D_j}\,dm_N\\
&=-N^2(N^2-1)c_N\,d(A,T_a,T_b)d(A,D_i,D_j)\\
&=-\frac{N^3}{2(N^2-4)}
       d(A,T_a,T_b)d(A,D_i,D_j),
\end{align*}
where the last line uses
$c_N=N/[2(N^2-1)(N^2-4)]$.
By definition, $d(A,T_a,T_b)=\langle\mathsf D_AT_a,T_b\rangle$
and similarly on the Cartan factor.  The entries therefore
factor as tensor products, yielding the real symmetric matrix
\begin{equation}\label{eq:flag-tensor-factor}
M_{\eta_N}=-\frac{N^3}{2(N^2-4)}
\mathsf D_A\otimes (\mathsf D_A|_{\mathfrak t_N}),
\end{equation}
where $\mathsf D_A$ is the self-adjoint operator on
$\mathfrak a_N$ uniquely determined by
\[
\langle\mathsf D_A B,C\rangle=d(A,B,C),
\qquad\text{equivalently}\qquad
\mathsf D_A B=\Proj_{\mathfrak a_N}(AB+BA).
\]
Here $\mathsf D_A|_{\mathfrak t_N}$ denotes the restriction of
this same operator to diagonal traceless matrices.  The factor
$N^3/[2(N^2-4)]$ follows from
$(2N)^2/4=N^2$,
$\dim\mathfrak a_N=N^2-1$, and the value of $c_N$ in
\eqref{eq:flag-cubic}.

We diagonalize both factors without leaving any undetermined rank.
Write $\mathbb C^N=\mathbb C^{N-1}\oplus\mathbb C$.
The following spaces are pairwise orthogonal and invariant under
$\mathsf D_A$:
\begin{enumerate}
\item traceless Hermitian matrices in the upper-left
$(N-1)\times(N-1)$ block, of real dimension
$(N-1)^2-1=N(N-2)$, on which $\mathsf D_A=2I$,
\item Hermitian off-diagonal blocks joining the first $N-1$
coordinates to the last, of real dimension $2(N-1)$,
on which $\mathsf D_A=-(N-2)I$,
\item the line spanned by $A$, on which
$\mathsf D_A=-2(N-2)I$.
\end{enumerate}
Indeed, for $B=\left(\begin{smallmatrix}C&z\\z^*&b\end{smallmatrix}\right)$
with $\operatorname{tr}C+b=0$,
\[
AB+BA=\begin{pmatrix}2C&-(N-2)z\\
-(N-2)z^*&-2(N-1)b\end{pmatrix}.
\]
Subtracting its scalar trace part gives the first two eigenvalues
on the indicated blocks and
$\mathsf D_AA=2\Proj_{\mathfrak a_N}(A^2)=-2(N-2)A$.
The dimensions sum to $N^2-1$, so
\begin{equation}\label{eq:flag-left-spectrum}
\operatorname{spec}(\mathsf D_A)=
\{2^{[N(N-2)]},\;[-(N-2)]^{[2(N-1)]},\;
[-2(N-2)]^{[1]}\}.
\end{equation}
(The superscripts in square brackets indicate multiplicities.)
On the Cartan subspace $\mathfrak t_N$, the first class restricts to
upper-left traceless diagonal matrices of dimension $N-2$,
and the third class remains a line. Thus
\begin{equation}\label{eq:flag-right-spectrum}
\operatorname{spec}(\mathsf D_A|_{\mathfrak t_N})
=\{2^{[N-2]},\;[-2(N-2)]^{[1]}\}.
\end{equation}
Neither factor has a zero eigenvalue. Their sign multiplicities are
\[
\begin{array}{c|cc}
&\dim E_+&\dim E_-\\ \hline
\mathsf D_A&N(N-2)&2N-1\\
\mathsf D_A|_{\mathfrak t_N}&N-2&1
\end{array}.
\]
Since the scalar in \eqref{eq:flag-tensor-factor} is negative,
equal-sign pairs give the negative eigenvalues of $M_{\eta_N}$
and opposite-sign pairs give its positive eigenvalues. Hence
\[
n_-(M_{\eta_N})=N(N-2)^2+2N-1,\qquad
n_+(M_{\eta_N})=(N-2)(3N-1).
\]
These counts sum to $(N^2-1)(N-1)$. The first count contains the
product of a quadratic-dimensional positive block in the adjoint
factor with a linear-dimensional positive block in the Cartan
multiplicity factor. The opposite-sign count is only quadratic.
This is the representation-theoretic source of the asymmetry.
For example, the two indices are $23,22$ at $N=4$ and $659,232$
at $N=10$.

Finally, the general local index theorem
(Theorem~\ref{thm:GH-index}) applies: on the orthogonal complement
of $E_1$ its defect form is uniformly coercive for sufficiently small
$|\varepsilon|$, and the reduced Schur complement is
$\varepsilon M_{\eta_N}+O(\varepsilon^2)$ on the fixed
finite-dimensional space $E_1$. Because $M_{\eta_N}$ is
invertible, this error cannot change its inertia for small
nonzero $\varepsilon$. Hence the full entropy-Hessian defect has
exactly the indices in \eqref{eq:flag-index}.
\end{proof}

The proof also establishes that the restriction of the first-order
defect to $E_1$ is nondegenerate. Its full signature is determined
by the two invariant cubic operators in \eqref{eq:flag-tensor-factor},
whose spectra are given in \eqref{eq:flag-left-spectrum} and
\eqref{eq:flag-right-spectrum}. In particular, the local index
changes under reversal of the perturbation for every $N\ge4$,
Corollary~\ref{cor:flag-asymmetry} makes this comparison explicit.

\begin{corollary}\label{cor:flag-su3}
For $SU(3)/T^2$, the first eigenspace has dimension $16$ and
multiplicity two for the adjoint $SU(3)$-type. For
$\eta_3(gT)=9|g_{33}|^2-3$, the first-order defect matrix has
eigenvalues $\pm54/5$ and $\pm27/5$, each with multiplicity four.
In particular, for both signs of sufficiently small
$\varepsilon\ne0$,
$\operatorname{ind}_-(\mathfrak d_{1+\varepsilon\eta_3})=8$.
\end{corollary}
\begin{proof}
For $N=3$, equations \eqref{eq:flag-left-spectrum} and
\eqref{eq:flag-right-spectrum} give eigenvalues
$\{2^{[3]},-1^{[4]},-2^{[1]}\}$ and
$\{2^{[1]},-2^{[1]}\}$ respectively. The constant in
\eqref{eq:flag-tensor-factor} is $-27/10$.
Multiplying these spectra gives the displayed values and the index.
\end{proof}

\begin{corollary}\label{cor:flag-asymmetry}
Let $E_{1,N}$ be the first positive Laplace eigenspace of
$\Fl_N$, and write 
\[
d_N=\dim E_{1,N}=(N^2-1)(N-1).
\]
For $0<t<1/[N(N-1)]$, the restrictions of the two entropy-Hessian defects
to $E_{1,N}$ are nondegenerate and satisfy
\begin{equation}\label{eq:flag-reversed-defect}
\left.\mathfrak d_{1-t\eta_N}\right|_{E_{1,N}}
=-\left.\mathfrak d_{1+t\eta_N}\right|_{E_{1,N}}.
\end{equation}
Their negative indices $I_N^+$ and $I_N^-$ are therefore complementary:
\begin{equation}\label{eq:flag-complementarity}
I_N^++I_N^-=d_N.
\end{equation}
For each fixed $N$, the same numbers are the negative indices
of the \emph{full} defect forms when $t>0$ is sufficiently small.
Moreover,
\begin{equation}\label{eq:flag-high-rank-rates}
\frac{I_N^-}{d_N}=\frac{3}{N}+O(N^{-2}),
\qquad
\frac{I_N^+}{d_N}=1-\frac{3}{N}+O(N^{-2})
\quad (N\to\infty).
\end{equation}
In particular, $I_N^\pm$ is the maximal dimension of a subspace
$V\subset E_{1,N}$ on which
\[
\mathfrak h_{1\pm t\eta_N}(u,u)<2N\,\mathfrak g_{1\pm t\eta_N}(u,u)
\]
for all $0\ne u\in V$.  The dimension imbalance does not measure
the magnitude of that inequality's failure. The parameter range
guaranteeing the \emph{full}-index identities may depend on $N$.
\end{corollary}
\begin{proof}
On $E_{1,N}$ the defect is exactly $\pm tM_{\eta_N}$ by
\eqref{eq:flag-first-defect}. The matrix is invertible, so reversing
its sign gives complementary indices. Theorem~\ref{thm:flag-index}
identifies these with the full indices for each fixed $N$ and
sufficiently small $t$. Finally,
\[
\frac{I_N^-}{d_N}
=\frac{3N^2-7N+2}{N^3-N^2-N+1}
=\frac3N+O(N^{-2}),
\]
and $I_N^+/d_N=1-I_N^-/d_N$.
\end{proof}

The asymmetry describes the number of directions below the
threshold for this family and this normal metric. It does not
classify all perturbations or invariant metrics on complete flags.
The local parameter range need not be uniform in $N$.

The same diagonalization also gives the first-order sizes of the
generalized Hessian eigenvalue shifts. Since the homogeneous Otto
metric on $E_1$ is $2NI$, the branch corresponding to eigenvalues
$s$ of $\mathsf D_A$ and $t$ of its Cartan restriction satisfies
\[
\kappa_{s,t}(\varepsilon)
=2N-\varepsilon\frac{N^2}{4(N^2-4)}st+O(\varepsilon^2).
\]
Thus the tensor-product computation yields both the exact
multiplicity of the crossings and their slopes.

\section{The circle: exact Fourier couplings}
\label{sec:circle-example}

\begin{example} Let $G=S^1$, $H=\{1\}$, use the metric $d\theta^2$, and normalize $dm=d\theta/(2\pi)$. The characters
\begin{equation}\label{eq:circle-modes}
\Phi_k(\theta)=e^{ik\theta},\qquad k\in\mathbb Z,
\end{equation}
are orthonormal and satisfy $-\Delta\Phi_k=k^2\Phi_k$. Therefore
\[
\|\nabla\Phi_k\|_{W,m}^2=k^2,
\qquad
\Hess^W\Ent_m(\nabla\Phi_k,\nabla\Phi_k)=k^4,
\]
and the infinitesimal Hessian rate is $k^2$.

Take
\begin{equation}\label{eq:circle-density}
\rho(\theta)=1+2\varepsilon\cos(q\theta)
=1+\varepsilon\Phi_q+\varepsilon\Phi_{-q},
\qquad 0<|\varepsilon|<\frac12,\quad q\geq1.
\end{equation}
The triple coefficient is
\begin{equation}\label{eq:circle-triple}
\mathcal C_{k\ell s}
=\int_0^{2\pi}e^{ik\theta}e^{-i\ell\theta}e^{is\theta}
\frac{d\theta}{2\pi}=\delta_{k-\ell+s,0}.
\end{equation}
Since
$\tfrac12(k^2+\ell^2-(\ell-k)^2)=k\ell$, the exact formulas become
\begin{align}
\mathfrak g_\rho(\Phi_k,\Phi_\ell)
&=k^2\delta_{k\ell}
+\varepsilon k\ell
(\delta_{\ell-k,q}+\delta_{\ell-k,-q}),
\label{eq:circle-G}\\
\mathfrak h_\rho(\Phi_k,\Phi_\ell)
&=k^4\delta_{k\ell}
+\varepsilon k^2\ell^2
(\delta_{\ell-k,q}+\delta_{\ell-k,-q}).
\label{eq:circle-H}
\end{align}
Thus, among nonconstant Fourier modes, the density couples exactly
those whose indices differ by $q$. Constant modes have zero gradient
and are excluded from these differential-form couplings.
For $q=2$,
\[
\mathfrak g_\rho(\Phi_1,\Phi_3)=3\varepsilon,
\qquad
\mathfrak h_\rho(\Phi_1,\Phi_3)=9\varepsilon,
\]
whereas both entries vanish for $(1,4)$. Under the left regular action
$(z\cdot f)(x)=f(z^{-1}x)$, the function $\Phi_k$ has type
$\chi_{-k}$. Hence
$\chi_{-k}\otimes\chi_{-\ell}^*=\chi_{\ell-k}$, while the dual of the
type of $\Phi_s$ is $\chi_s$.  The selection rule is therefore
$s=\ell-k$, or $k-\ell+s=0$, exactly the condition in
\eqref{eq:circle-triple}.
\end{example}

The entries of the Otto metric above are the complex Fourier version of established Fourier-frame formulas.  In particular, Proposition 3.1 of \cite{GomesRodriguesSanMartin2026} gives explicit coefficients at arbitrary smooth circle densities in a real sine--cosine frame. The purpose of the present example is to test the \emph{entropy-Hessian and negative-index} formulas, not to reintroduce the Fourier metric as a new result.
\begin{corollary}\label{cor:circle-obstruction}
For the density \eqref{eq:circle-density}, every admissible
$\varepsilon\neq0$ yields a negative direction of
$\mathfrak d_\rho$.  For sufficiently small nonzero $\varepsilon$,
\begin{equation}\label{eq:circle-index}
\operatorname{ind}_-(\mathfrak d_\rho)
=\begin{cases}1,&q=2,\\2,&q\neq2.\end{cases}
\end{equation}
\end{corollary}
\begin{proof}
The complex first eigenspace is
$E_1=\operatorname{span}\{\Phi_1,\Phi_{-1}\}$.
For $q=2$, \eqref{eq:circle-G}--\eqref{eq:circle-H} yield
\[
M_\eta=\begin{pmatrix}0&2\\2&0\end{pmatrix}
\]
for $\eta=2\cos(2\theta)$, so $M_\eta$ is invertible and
$n_-(\varepsilon M_\eta)=1$.  Theorem~\ref{thm:GH-index}
gives the first value in \eqref{eq:circle-index}.
If $q\neq2$, then $M_\eta=0$.  The mixed coupling of
$\Phi_{1}$ with $\Phi_{q+1}$ equals
$\varepsilon q(q+1)$, and the conjugate coupling of
$\Phi_{-1}$ with $\Phi_{-q-1}$ has the same nonzero magnitude.
These two outputs are orthogonal and distinct, so $B_\eta$ has
rank two on $E_1$.  Theorem~\ref{thm:GH-index} gives index two.
For either value of $q$, the asserted negative direction for every
admissible parameter follows already from the nonzero first-row
pairing and Lemma~\ref{lem:tracefree} or the two-dimensional
mixed determinant argument.
\end{proof}

\appendix
\section{An explicit \texorpdfstring{$SU(2)$}{SU(2)} coupling calculation}
\label{sec:su2-calculation}

The channel calculation in Appendix~\ref{sec:su2-detection} uses
normalized Clebsch--Gordan coefficients, not just the abstract
fusion rule. We derive those coefficients directly, including the
normalization of the symmetric-power basis, the highest-weight
intertwiners, and the resulting Haar triple products.

\subsection*{Symmetric powers and the orthonormal weight basis}
Let $U=\mathbb C^2$ with orthonormal basis $z,w$. The irreducible
$SU(2)$-representation $V_k=\operatorname{Sym}^k(U)$ is realized by
homogeneous polynomials of degree $k$. Its dimension is $k+1$, and
the Peter--Weyl block of its matrix coefficients has dimension
$(k+1)^2$ by Schur orthogonality. The nonnegative operator $L=-\Delta$ acts on
this block by the eigenvalue
\begin{equation}\label{eq:SU2-spectrum}
L|_{\mathcal W_k}=\lambda_k I,\qquad \lambda_k=k(k+2).
\end{equation}
The tensor-product decomposition is
\begin{equation}\label{eq:SU2-CG}
V_n\otimes V_m=\bigoplus_{r=0}^{\min\{n,m\}}V_{n+m-2r}.
\end{equation}

For $0\le j\le k$, let $T_I$ be the elementary tensor having $w$
in the positions $I\subset\{1,\ldots,k\}$ and $z$ elsewhere.
Distinct $T_I$ are orthonormal. The symmetrization of
$z^{\otimes(k-j)}\otimes w^{\otimes j}$ is therefore
\[
S_{k,j}=\frac1{k!}\sum_{\sigma\in\mathfrak S_k}
P_\sigma(z^{\otimes(k-j)}\otimes w^{\otimes j})
=\frac1{\binom{k}{j}}\sum_{|I|=j}T_I,
\]
where each $T_I$ occurs $j!(k-j)!$ times. It follows that
\[
\|S_{k,j}\|^2
=\frac1{\binom{k}{j}^2}\sum_{|I|=j}\|T_I\|^2
=\frac1{\binom{k}{j}}.
\]
The weight spaces are orthogonal, hence the normalized basis is
\begin{equation}\label{eq:su2-orthonormal-basis}
e_j^{(k)}=\sqrt{\binom{k}{j}}\,z^{k-j}w^j,
\qquad j=0,\ldots,k.
\end{equation}
Here the polynomial notation uses the invariant inner product
induced by the symmetric tensors, so the normalization in
\eqref{eq:su2-orthonormal-basis} is part of the calculation.

Let $H_0,E,F$ be the usual generators of
$\mathfrak{sl}_2(\mathbb C)$, acting on polynomials by
\[
H_0=z\partial_z-w\partial_w,\qquad
E=z\partial_w,\qquad F=w\partial_z.
\]
The Leibniz rule gives
\begin{align*}
H_0(z^{k-j}w^j)&=(k-2j)z^{k-j}w^j,\\
E(z^{k-j}w^j)&=jz^{k-j+1}w^{j-1},\\
F(z^{k-j}w^j)&=(k-j)z^{k-j-1}w^{j+1}.
\end{align*}
Using the binomial ratios
\[
\frac{\binom{k}{j}}{\binom{k}{j-1}}=\frac{k-j+1}{j},
\qquad
\frac{\binom{k}{j}}{\binom{k}{j+1}}=\frac{j+1}{k-j},
\]
we obtain the exact action on the normalized basis:
\begin{equation}\label{eq:su2-raising}
\begin{aligned}
H_0e_j^{(k)}&=(k-2j)e_j^{(k)},\\
Ee_j^{(k)}&=\sqrt{j(k-j+1)}\,e_{j-1}^{(k)},\\
Fe_j^{(k)}&=\sqrt{(j+1)(k-j)}\,e_{j+1}^{(k)}.
\end{aligned}
\end{equation}
Terms with out-of-range indices vanish. On $V_1\otimes V_k$,
each generator acts by $X\otimes I+I\otimes X$.

\subsection*{The two normalized intertwiners}
Set $k=n+1$. The Clebsch--Gordan decomposition gives
$V_1\otimes V_k=V_{k+1}\oplus V_{k-1}$. We construct the two
isometric embeddings
\[
J_+:V_{k+1}\longrightarrow V_1\otimes V_k,
\qquad
J_-:V_{k-1}\longrightarrow V_1\otimes V_k.
\]
The vector $e_0^{(1)}\otimes e_0^{(k)}$ is the unit
highest-weight vector of weight $k+1$. Fix the phase by setting
\[
J_+e_0^{(k+1)}=e_0^{(1)}\otimes e_0^{(k)}.
\]
At weight $k+1-2j$, write
\[
J_+e_j^{(k+1)}
=a_j e_0^{(1)}\otimes e_j^{(k)}
+b_j e_1^{(1)}\otimes e_{j-1}^{(k)}.
\]
The relation $FJ_+=J_+F$, using \eqref{eq:su2-raising},
requires
\begin{align*}
a_j\sqrt{(j+1)(k-j)}
 &=\sqrt{(j+1)(k+1-j)}\,a_{j+1},\\
a_j+b_j\sqrt{j(k-j+1)}
 &=\sqrt{(j+1)(k+1-j)}\,b_{j+1}.
\end{align*}
The initial data are $a_0=1$ and $b_0=0$. The first recurrence
telescopes to
\[
a_j=\prod_{r=0}^{j-1}\sqrt{\frac{k-r}{k+1-r}}
=\sqrt{\frac{k+1-j}{k+1}}.
\]
Inserting this into the second recurrence gives, inductively,
$b_j=\sqrt{j/(k+1)}$. Indeed its left-hand side becomes
\[
\sqrt{\frac{k+1-j}{k+1}}
+\sqrt{\frac{j}{k+1}}\sqrt{j(k-j+1)}
=(j+1)\sqrt{\frac{k+1-j}{k+1}},
\]
which equals the right-hand side with
$b_{j+1}=\sqrt{(j+1)/(k+1)}$. Thus
\begin{equation}\label{eq:su2-Jplus}
\begin{aligned}
J_+e_j^{(k+1)}={}&\sqrt{\frac{k+1-j}{k+1}}\,
 e_0^{(1)}\otimes e_j^{(k)}\\
&+\sqrt{\frac{j}{k+1}}\,
 e_1^{(1)}\otimes e_{j-1}^{(k)}.
\end{aligned}
\end{equation}

For the lower summand, consider the weight-$k-1$ vector
\[
v=a e_0^{(1)}\otimes e_1^{(k)}
+b e_1^{(1)}\otimes e_0^{(k)}.
\]
The equation $Ev=0$ gives $a\sqrt{k}+b=0$, choosing $a>0$
and $\|v\|=1$ yields
\[
J_-e_0^{(k-1)}
=\frac1{\sqrt{k+1}}e_0^{(1)}\otimes e_1^{(k)}
-\sqrt{\frac{k}{k+1}}e_1^{(1)}\otimes e_0^{(k)}.
\]
Write the remaining weight vectors as
\[
J_-e_j^{(k-1)}
=c_j e_0^{(1)}\otimes e_{j+1}^{(k)}
-d_j e_1^{(1)}\otimes e_j^{(k)}.
\]
The intertwining relation $FJ_-=J_-F$ gives the two recurrences
\begin{align*}
c_j\sqrt{(j+2)(k-j-1)}
 &=\sqrt{(j+1)(k-1-j)}\,c_{j+1},\\
c_j-d_j\sqrt{(j+1)(k-j)}
 &=-\sqrt{(j+1)(k-1-j)}\,d_{j+1}.
\end{align*}
Starting from $c_0=1/\sqrt{k+1}$ and
$d_0=\sqrt{k/(k+1)}$, these yield
\[
c_j=\sqrt{\frac{j+1}{k+1}},\qquad
 d_j=\sqrt{\frac{k-j}{k+1}}.
\]
For example, substituting these expressions in the second
recurrence reduces its two sides to
\[
-(k-1-j)\sqrt{\frac{j+1}{k+1}}
=-\sqrt{(j+1)(k-1-j)}
\sqrt{\frac{k-1-j}{k+1}}.
\]
Hence
\begin{equation}\label{eq:su2-Jminus}
\begin{aligned}
J_-e_j^{(k-1)}={}&\sqrt{\frac{j+1}{k+1}}\,
 e_0^{(1)}\otimes e_{j+1}^{(k)}\\
&-\sqrt{\frac{k-j}{k+1}}\,
 e_1^{(1)}\otimes e_j^{(k)}.
\end{aligned}
\end{equation}
The ranges are $0\le j\le k+1$ in \eqref{eq:su2-Jplus} and
$0\le j\le k-1$ in \eqref{eq:su2-Jminus}. In each formula,
the squared coefficients sum to one, and different weight vectors
are orthogonal. At a common weight, the inner product between
$J_+e_{j+1}^{(k+1)}$ and $J_-e_j^{(k-1)}$ is
\[
\frac{\sqrt{(j+1)(k-j)}}{k+1}
-\frac{\sqrt{(j+1)(k-j)}}{k+1}=0.
\]
The relations with $H_0$ and $F$ have been verified.
Substitution of the same coefficients into the raising formulas in
\eqref{eq:su2-raising} also gives $EJ_\pm=J_\pm E$.
Thus both maps are $\mathfrak{sl}_2(\mathbb C)$-intertwiners. Finally,
\[
\dim V_{k+1}+\dim V_{k-1}
=(k+2)+k=2(k+1)=\dim(V_1\otimes V_k),
\]
so the two images exhaust the tensor product. These are the
normalized intertwining coefficients rather than only the abstract
fusion multiplicities.

\subsection*{Haar triple products and the defect factors}
For a fixed orthonormal basis of $V_r$, write
\[
\Phi_{r,i,a}(g)
=\sqrt{r+1}\,\langle\pi_r(g)e_a^{(r)},e_i^{(r)}\rangle.
\]
For $\ell=k+1$ or $k-1$, write $J_\ell=J_+$ or $J_-$,
respectively, and set
\[
C^{(\ell)}_{ir,j}
=\langle J_\ell e_j^{(\ell)},
 e_i^{(1)}\otimes e_r^{(k)}\rangle.
\]
Schur orthogonality and the intertwiner projection give
\begin{equation}\label{eq:su2-CG-contraction}
\int_{SU(2)}\Phi_{1,i,a}\Phi_{k,r,c}
\overline{\Phi_{\ell,j,b}}\,dg
=\sqrt{\frac{2(k+1)}{\ell+1}}
\overline{C^{(\ell)}_{ac,b}}\,C^{(\ell)}_{ir,j},
\qquad \ell=k\pm1.
\end{equation}
This is the one-copy specialization of
\eqref{eq:CG-coordinate-contraction}. For $k=2$, the upper
highest-weight identity gives
\begin{equation}\label{eq:su2-example-upper-CG}
\begin{aligned}
\int_{SU(2)}\Phi_{1,0,0}\Phi_{2,0,0}
\overline{\Phi_{3,0,0}}\,dg
&=\sqrt{\frac{2\cdot3}{4}}\cdot1\cdot1
=\sqrt{\frac32}.
\end{aligned}
\end{equation}
For the lower channel, \eqref{eq:su2-Jminus} yields
\[
J_-e_0^{(1)}=\frac{e_0^{(1)}\otimes e_1^{(2)}}{\sqrt3}
-\sqrt{\frac23}\,e_1^{(1)}\otimes e_0^{(2)}.
\]
Thus the other explicit triple product is
\begin{equation}\label{eq:su2-example-lower-CG}
\int_{SU(2)}\Phi_{1,0,0}\Phi_{2,1,1}
\overline{\Phi_{1,0,0}}\,dg
=\sqrt{\frac{2\cdot3}{2}}\left(\frac1{\sqrt3}\right)^2
=\frac1{\sqrt3}.
\end{equation}
Both numbers can also be checked directly from Haar integration.
For Haar-distributed $g\in SU(2)$, the variable $X=|g_{11}|^2$
is uniform on $[0,1]$, the relevant matrix coefficients are
$g_{11}$, $g_{11}^2$, $g_{11}^3$, and $2X-1$. Therefore
\begin{align*}
\int_{SU(2)}\Phi_{1,0,0}\Phi_{2,0,0}
\overline{\Phi_{3,0,0}}\,dg
 &=2\sqrt6\int_0^1X^3\,dX=\sqrt{\frac32},\\
\int_{SU(2)}\Phi_{1,0,0}\Phi_{2,1,1}
\overline{\Phi_{1,0,0}}\,dg
 &=2\sqrt3\int_0^1X(2X-1)\,dX=\frac1{\sqrt3}.
\end{align*}
Finally the entropy-Hessian defect multiplies a permitted triple
product by the Laplace spectral factor from
Lemma~\ref{lem:first-row}. For $u\in\mathcal W_1$,
$v\in\mathcal W_{n+1}$, and a density component
$\eta_r\in\mathcal W_r$, that factor equals
\[
\frac{(\lambda_{n+1}-\lambda_r)^2-\lambda_1^2}{4}.
\]
Since $\lambda_j=j(j+2)$, the two permitted density modes are
$r=n$ and $r=n+2$, with explicit values
\begin{align*}
r=n:\quad
\frac{(2n+3)^2-9}{4}&=n(n+3),\\
r=n+2:\quad
\frac{(2n+5)^2-9}{4}&=(n+1)(n+4).
\end{align*}
Combining these numbers with the two Haar contractions yields the
channel formula in Lemma~\ref{lem:SU2-two-channels}. The
calculation therefore distinguishes which types are allowed, which
actual Clebsch--Gordan entries are nonzero, and the numerical
spectral factors multiplying those entries.

\section{Representation checks on \texorpdfstring{$SU(2)$}{SU(2)}}
\label{sec:su2-detection}

Equip $SU(2)$ with the bi-invariant metric of
\eqref{eq:SU2-spectrum}, identifying it with the unit round $S^3$.
Its full Peter--Weyl component $\mathcal W_n$ has dimension
$(n+1)^2$ and coincides with the complexified spherical harmonics
of degree $n$. The real characters $\chi_n=\operatorname{tr}\pi_n$
are orthonormal among class functions and satisfy
\[
\chi_n\chi_m=\sum_{j=0}^{\min(n,m)}\chi_{n+m-2j}.
\]
The calculations below check the spherical results in
Clebsch--Gordan coordinates and give an exact negative test
potential.

\begin{example}\label{ex:su2-fusion}
Set
\begin{equation}\label{eq:SU2-density}
\rho=1+\varepsilon\chi_2,\qquad
-\frac13<\varepsilon<1,\quad\varepsilon\ne0.
\end{equation}
Since $\chi_2(\theta)=1+2\cos(2\theta)\in[-1,3]$, the displayed
interval is exactly the positivity range, and $\rho$ is a smooth
positive probability density. Since
$V_1\otimes V_3=V_4\oplus V_2$, multiplication of characters gives
$\chi_1\chi_3=\chi_4+\chi_2$. Therefore character orthogonality yields
\[
\mathcal C_{1,3,2}
=\int_{SU(2)}\chi_1\chi_3\chi_2\,dg
=\langle\chi_4+\chi_2,\chi_2\rangle_{L^2}=1.
\]
With
$\lambda_1=3$, $\lambda_2=8$, and $\lambda_3=15$,
\begin{align}
\mathfrak g_\rho(\chi_1,\chi_3)
&=\frac{\varepsilon}{2}(3+15-8)=5\varepsilon,
\label{eq:SU2-nonzero-G}\\
\mathfrak h_\rho(\chi_1,\chi_3)
&=\frac{\varepsilon}{4}(3+15-8)^2=25\varepsilon.
\label{eq:SU2-nonzero-H}
\end{align}
Thus the allowed $V_2$ channel produces an explicitly nonzero entry.
\end{example}

\paragraph{Parity and exact decoupling.}
Let $Z=\{\pm I\}$ be the center of $SU(2)$ and suppose the smooth
positive density $\rho$ is $Z$-invariant.  Since the center acts on
$V_\ell$ by $(-1)^\ell$, its Peter--Weyl expansion contains only
even degrees $\ell$.  The fusion rule \eqref{eq:SU2-CG} says that
every type in $V_n\otimes V_m$ has the parity of $n+m$.
If $n-m$ is odd, all the types in the product are odd and thus none
can pair with a density coefficient.  The constant-density terms also
vanish because $n\neq m$.  Hence the exact formulas
\eqref{eq:exact-G}--\eqref{eq:exact-H} imply
\begin{equation}\label{eq:SU2-all-density-vanishing}
\mathfrak g_\rho(\mathcal W_n,\mathcal W_m)=0,
\qquad \mathfrak h_\rho(\mathcal W_n,\mathcal W_m)=0
\quad\text{if $n-m$ is odd}.
\end{equation}
For instance, $\mathcal W_1$ and $\mathcal W_2$ remain uncoupled,
whereas $\mathcal W_1$ and $\mathcal W_3$ may couple, as in
\eqref{eq:SU2-nonzero-G}--\eqref{eq:SU2-nonzero-H}.
This is a direct application of representation selection, rather
than an additional theorem about quotient integration.

\begin{lemma}\label{lem:SU2-two-channels}
For $n\geq1$, $u\in\mathcal W_1$, and $v\in\mathcal W_{n+1}$,
write $\rho=1+\sum_{k\geq1}\eta_k$, $\eta_k\in\mathcal W_k$.
Then
\begin{equation}\label{eq:SU2-Dchannels}
\begin{aligned}
\mathfrak d_\rho(u,v)
={}&n(n+3)\int_{SU(2)}\eta_nu\overline v\,dm+(n+1)(n+4)\int_{SU(2)}\eta_{n+2}u\overline v\,dm.
\end{aligned}
\end{equation}
Moreover,
\begin{equation}\label{eq:SU2-noncancel}
\mathfrak d_\rho(\mathcal W_1,\mathcal W_{n+1})=0
\quad\Longleftrightarrow\quad \eta_n=\eta_{n+2}=0.
\end{equation}
\end{lemma}
\begin{proof}
The product $V_1\otimes V_{n+1}=V_n\oplus V_{n+2}$ permits
precisely these two density components. Since
$\lambda_k=k(k+2)$, the numerical factors in
Lemma~\ref{lem:first-row} are $n(n+3)$ and $(n+1)(n+4)$.
For noncancellation use the left--right action of
$SU(2)\times SU(2)$, under which
$\mathcal W_k\simeq V_k^*\boxtimes V_k$ is irreducible.
Each map
\[
\eta_k\longmapsto
\left((u,v)\longmapsto\int\eta_ku\overline v\,dm\right),
\qquad k=n,n+2,
\]
is equivariant and nonzero, as
$\int\chi_k\chi_1\chi_{n+1}\,dm=1$.
Its kernel is therefore zero. The two images are inequivalent
irreducible submodules, so their intersection is zero.
This proves \eqref{eq:SU2-noncancel}, independently verifying
Lemma~\ref{lem:sphere-channels} for $d=3$.
\end{proof}

\begin{corollary}\label{cor:SU2-index}
For any smooth positive probability density on $SU(2)$,
\begin{equation}\label{eq:SU2-rigidity}
\mathfrak h_\rho\geq3\mathfrak g_\rho
\quad\Longleftrightarrow\quad \rho\equiv1.
\end{equation}
For any nonzero smooth real mean-zero $\eta$, write its
degree-two component as $x^{\mathsf T}Q_\eta x$ on $S^3$.
In the orthonormal basis $u_i=2x_i$, $0\leq i\leq3$,
\begin{equation}\label{eq:SU2-obstruction-matrix}
M_\eta=\frac43 Q_\eta.
\end{equation}
For sufficiently small nonzero $\varepsilon$,
\begin{equation}\label{eq:SU2-degree2-index}
\operatorname{ind}_-(\mathfrak d_{1+\varepsilon\eta})
=n_-(\varepsilon Q_\eta)+\dim\ker Q_\eta,\qquad
\operatorname{nul}(\mathfrak d_{1+\varepsilon\eta})=0.
\end{equation}
In particular, for $0\ne\eta\in\mathcal W_n^{\mathbb R}$ with
$n\geq1$ and $n\ne2$, the local index is four.
\end{corollary}
\begin{proof}
These assertions are Theorem~\ref{thm:sphere-index} and
\eqref{eq:sphere-moment-matrix} in dimension three.
\end{proof}

For the density in Example~\ref{ex:su2-fusion},
$\chi_2=4x_0^2-1$ gives
$Q_{\chi_2}=\operatorname{diag}(3,-1,-1,-1)$.
Thus the full local index is three for positive $\varepsilon$
and one for negative $\varepsilon$. The mixed character pairing
in that example is consistent with this first-order splitting.

\begin{example}
\label{ex:SU2-exact-quadratic}
For $\eta=\chi_3$ and $\rho_\varepsilon=1+\varepsilon\chi_3$,
consider the real central potential
\[
u_\varepsilon=\chi_1-\varepsilon
\left(\tfrac14\chi_2+\tfrac1{28}\chi_4\right).
\]
The fusion relation $\chi_3\chi_1=\chi_2+\chi_4$ and the
spectral factors from \eqref{eq:resonance-coupling} give
\[
\mathfrak d_\eta(\chi_1,\chi_2)=10,
\qquad \mathfrak d_\eta(\chi_1,\chi_4)=18.
\]
The homogeneous diagonal values are
$\mathfrak d_1(\chi_2,\chi_2)=8(8-3)=40$ and
$\mathfrak d_1(\chi_4,\chi_4)=24(24-3)=504$.
The $\eta$-weighted terms between the even-degree characters
$\chi_2,\chi_4$ vanish by parity, as does
$\mathfrak d_\eta(\chi_1,\chi_1)$.  Therefore the full quadratic
form can be evaluated, without an asymptotic remainder:
\begin{align*}
\mathfrak d_{\rho_\varepsilon}(u_\varepsilon,u_\varepsilon)
&=\varepsilon^2\left(\frac{40}{16}+\frac{504}{28^2}
-2\left(\frac{10}{4}+\frac{18}{28}\right)\right)\\
&=-\frac{22}{7}\varepsilon^2.
\end{align*}
Thus two permitted Clebsch--Gordan channels yield an explicitly
negative Wasserstein Hessian defect for every admissible nonzero
$\varepsilon$, not merely a formal perturbation series.
\end{example}

\section*{Acknowledgements}

This work was financed in part by the Coordena\c{c}\~{a}o de Aperfei\c{c}oamento de
Pessoal de N\'{i}vel Superior - Brasil (CAPES) - Finance Code 001. C. S. R. has been partially supported by S\~{a}o Paulo Research Foundation (FAPESP): grant \#2018/13481-0, and grant \#2020/04426-6. The opinions, hypotheses and conclusions or recommendations expressed in this work are the responsibility of the authors and do not necessarily reflect the views of FAPESP.

G. C. B. and C. S. R. would like to acknowledge support from the Max Planck Society, Germany, through the award of a Max Planck Partner Group for Geometry and Probability in Dynamical Systems.


\begin{thebibliography}{99}
\bibitem{FultonHarris1991}
W.~Fulton and J.~Harris,
\emph{Representation Theory: A First Course},
Graduate Texts in Mathematics 129, Springer, 1991.

\bibitem{Procesi1976}
C.~Procesi,
\emph{The invariant theory of $n\times n$ matrices},
Adv. Math. \textbf{19} (1976), no.~3, 306--381,
doi:10.1016/0001-8708(76)90027-X.

\bibitem{Yamaguchi1979}
S.~Yamaguchi,
\emph{Spectra of flag manifolds},
Mem. Fac. Sci. Kyushu Univ. Ser. A \textbf{33} (1979), no.~1, 95--112,
 doi:10.2206/kyushumfs.33.95.

\bibitem{Haber2021}
H.~E.~Haber,
\emph{Useful relations among the generators in the defining and
adjoint representations of $SU(N)$},
SciPost Phys. Lect. Notes \textbf{21} (2021),
doi:10.21468/SciPostPhysLectNotes.21.


\bibitem{GomesRodriguesInvariant2026}
A.~M.~de S\'a Gomes and C.~S.~Rodrigues,
\emph{Displacement convexity of invariant measures and curvature bounds
of isometric actions}, arXiv:2310.15332 (revised 2026).

\bibitem{GomesRodriguesSanMartin2026}
A.~M.~de S\'a Gomes, C.~S.~Rodrigues, and L.~A.~B.~San Martin,
\emph{The Wasserstein space of the circle is not flat},
arXiv:2610.00625 (2026).

\bibitem{Farashahi2017}
A.~Ghaani Farashahi,
\emph{Peter--Weyl theorem for homogeneous spaces of compact groups},
Int. J. Anal. Appl. \textbf{13} (2017), no.~1, 22--31.

\bibitem{Bakry1991}
D.~Bakry,
\emph{In\'egalit\'es de Sobolev faibles : un crit\`ere $\Gamma_2$},
S\'eminaire de Probabilit\'es \textbf{25} (1991), 234--261.

\bibitem{Kato1995}
T.~Kato,
\emph{Perturbation Theory for Linear Operators},
Classics in Mathematics, Springer, 1995.

\bibitem{Maddocks1988}
J.~H.~Maddocks,
\emph{Restricted quadratic forms, inertia theorems, and the Schur complement},
Linear Algebra Appl. \textbf{108} (1988), 1--36.

\bibitem{AxlerBourdonRamey2001}
S.~Axler, P.~Bourdon, and W.~Ramey,
\emph{Harmonic Function Theory}, 2nd ed.,
Graduate Texts in Mathematics 137, Springer, 2001.

\bibitem{BrockerDieck1985}
T.~Br\"ocker and T.~tom Dieck,
\emph{Representations of Compact Lie Groups},
Graduate Texts in Mathematics 98, Springer, 1985.

\bibitem{Lott2008}
J.~Lott,
\emph{Some geometric calculations on Wasserstein space},
Comm. Math. Phys. \textbf{277} (2008), 423--437.

\bibitem{RenesseSturm2005}
M.-K.~von Renesse and K.-T.~Sturm,
\emph{Transport inequalities, gradient estimates, entropy and Ricci curvature},
Comm. Pure Appl. Math. \textbf{58} (2005), no.~7, 923--940.

\bibitem{LottVillani2009}
J.~Lott and C.~Villani,
\emph{Ricci curvature for metric-measure spaces via optimal transport},
Ann. of Math. (2) \textbf{169} (2009), 903--991.

\bibitem{Otto2001}
F.~Otto,
\emph{The geometry of dissipative evolution equations: the porous medium equation},
Comm. Partial Differential Equations \textbf{26} (2001), 101--174.

\bibitem{Villani2009}
C.~Villani,
\emph{Optimal Transport: Old and New},
Grundlehren der mathematischen Wissenschaften 338, Springer, 2009.
\end{thebibliography}
\end{document}